\documentclass[epj]{svjour}
\usepackage{graphics}
\usepackage{amssymb}
\usepackage{tikz}
\begin{document}
\title{Simulations of Kinetic Electrostatic Electron Nonlinear (KEEN) Waves with Variable Velocity Resolution Grids and High-Order Time-Splitting}
\titlerunning{KEEN waves simulations}
\authorrunning{B. Afeyan et al}

\author{B. Afeyan\inst{1},  F. Casas\inst{2}, N. Crouseilles\inst{3}, A. Dodhy\inst{4},
E. Faou\inst{3}, M. Mehrenberger\inst{4,5}, E. Sonnendr\"ucker\inst{4,6}
\thanks{ A part of this work was carried out using the HELIOS supercomputer system at Computational Simulation Centre of International Fusion Energy Research Centre (IFERC-CSC), Aomori, Japan, under the Broader Approach collaboration between Euratom and Japan, implemented by Fusion for Energy and JAEA.
}
\thanks{BA would like to acknowledge the financial assistance of DOE OFES HEDP program through a subcontract via UCSD}
\thanks{MM would like to thank specially Edwin Chacon-Golcher and also Matthieu Haefele and Pierre Navaro for their help
on the development of the code in the SELALIB library.}
\thanks{This work was partly supported by the ERC Starting Grant Project GEOPARDI.}
%
}                     
%
%
\institute{Polymath Research Inc., Pleasanton, CA 94566 
\and IMAC, Departament de Matem\`atiques, Universitat Jaume I, 12071-Castell\'on, Spain
\and Inria-Rennes (IPSO team) and IRMAR, Rennes, France
\and Max-Planck-Institut f\"ur Plasmaphysik, Garching, Germany  \and IRMA, Universit\'e de Strasbourg, France
\and Mathematics center, TU Munich,  Garching, Germany}
\date{Received: date / Revised version: date}
%
\abstract{
KEEN waves are non-stationary, nonlinear, self-organized asymptotic states in Vlasov plasmas. They lie outside the precepts of linear theory or perturbative analysis, unlike electron plasma waves or ion acoustic waves. Steady state, nonlinear constructs such as BGK modes also do not apply. The range in velocity that is strongly perturbed by KEEN waves depends on the amplitude and duration of the ponderomotive force generated by two crossing laser beams, for instance, used to drive them. Smaller amplitude drives manage to devolve into multiple highly-localized vorticlets, after the drive is turned off, and may eventually succeed to coalesce into KEEN waves. Fragmentation once the drive stops, and potential eventual remerger, is a hallmark of the weakly driven cases. A fully formed (more strongly driven) KEEN wave has one dominant vortical core. But it also involves fine scale complex dynamics due to shedding and merging of smaller vortical structures with the main one. Shedding and merging of vorticlets are involved in either case, but at different rates and with different relative importance. The narrow velocity range in which one must maintain sufficient resolution in the weakly driven cases, challenges fixed velocity grid numerical schemes. What is needed is the capability of resolving locally in velocity while maintaining a coarse grid outside the highly perturbed region of phase space. We here report on a new Semi-Lagrangian Vlasov-Poisson solver based on  conservative non-uniform cubic splines in velocity that tackles this problem head on. An additional feature of our approach is the use of a new high-order time-splitting scheme which allows much longer simulations per computational effort. This is needed for low amplitude runs. There, global coherent structures take a long time to set up, such as KEEN waves, if they do so at all. The new code's performance is compared to uniform grid simulations and the advantages are quantified. The birth pains associated with weakly driven KEEN waves are captured in these simulations. Canonical KEEN waves with ample drive are also treated using these advanced techniques. They will allow the efficient simulation of KEEN waves in multiple dimensions, which will be tackled next, as well as generalizations to Vlasov-Maxwell codes. These are essential for pursuing the impact of KEEN waves in high energy density plasmas and in inertial confinement fusion applications. More generally, one needs a fully-adaptive grid-in-phase-space method which could handle all small vorticlet dynamics whether pealing off or remerging. Such fully adaptive grids would have to be computed sparsely in order to be viable. This two-velocity grid method is a concrete and fruitful step in that direction. 
%
} 
\maketitle
\section{Introduction}
\label{intro}

Kinetic Electrostatic Electron Nonlinear (KEEN) Waves were discovered in 2002 \cite{bedros04}. The impetus came from the examination of the validity of claims made regarding electron acoustic waves, EAW, and their relationship to electron plasma waves (EPW) in the nonlinear, kinetic evolution of the Stimulated Raman Scattering (SRS) instability \cite{bedros_ref2,bedros_ref3}. EAW and similar stationary constructs are single mode prescriptions, relying on a cobbled electron distribution function. They are assumed to be far from a Maxwellian, with a zero slope somewhere in the active (wave perturbed) region in velocity. They also invariably  lie on some particular dispersion curve. None of these feature pertain to KEEN waves. They have no infinitesimal amplitude, nor a fluid limit, they are generically a self-consistent multi-mode response of the plasma. They do not lie on a dispersion curve but spatial and temporal harmonics do arise in pairs (and even triplets, if velocity harmonics are also counted). KEEN waves are phase-locked, multiple-harmonic phase-space structures which are non-stationary. These coherent wave structures accompany some deeply trapped, and other chaotic or trapping-untrapping and distantly re-trapping oscillations.  Since that initial discovery, much work has been done to decipher the physics of KEEN waves. 
KEEN waves are most readily driven by the ponderomotive force generated by the optical mixing of a pair of laser beams \cite{bedros04,bedros14b} in an otherwise Maxwellian plasma. These two crossing laser beams drive a wave at their difference frequency and wavenumber. By changing the frequencies and wave numbers of the two laser beams in a given density and temperature plasma, we can drive KEEN waves anywhere in the Brillouin $(\omega, k)$ plane. Linear wave analysis of the Vlasov-Poisson system of equations, shows that for a Maxwellian plasma, for instance, resonant waves can only live on dispersion curves, which are familiar from plasma physics textbooks and are closely related to Landau's original work on collisionless damping of EPWs \cite{landau}.

Outside these curves in the $(\omega, k)$ plane, no waves were deemed possible unless they were heavily damped (given that the dispersion relation for electron plasma oscillations is a transcendental equation with an infinite number of zeros, one of them being the least damped root). Thus, the linear theory perspective defined a spectral gap in plasma theory. On the other hand, nonlinear stationary states were known and referred to as  BGK modes \cite{bgk}. They are predicated on one's ability to cobble a distribution function in phase space that can accommodate such a mode (This is done by one of two ways. Either by guessing the modified distribution function form, such as flatness at the phase velocity of the wave, and then determining the self-consistent electric field.  Or, by limiting the electric field to some  reductive form, such as it being single mode, constant amplitude, stationarity, and solving for the modified distribution function. In that case, the distribution function would be composed of a temporally frozen set of trapped and untrapped particle orbits).  These stationary modes (in, at best, some Galilean frame of reference and constituting a function of the canonically conjugate variable to time, namely, energy) require strong and permanent distribution function modifications. Thus, they beg the question: how could they come about starting with a  Maxwellian equilibrium plasma? Or, whether some other nonstationary states would predominate much before these fixed points were ever reached. The discovery of KEEN waves \cite{bedros14b} answers these questions in favor of the latter. Since KEEN waves form anywhere in the $(\omega, k)$ plane with a sufficiently long and strong drive, and not just on, or in the vicinity of the EPW and EAW dispersion curves, KEEN waves remove the delicate nature of linear, resonant modes as the only locations where surviving nonlinear structures can be constructed via adiabatic but persistent perturbations. 

We now know that EAWs and BGK modes are not the only type of nonlinear structure to be expected in a time dependent setting. KEEN waves demonstrate that the requirement of stationarity is too restrictive and that constantly evolving and adapting phase space structures may be ubiquitous, nonlinear states of self-organization of a plasma. The basic physics is that self-consistent electric field structures made up of multiple  phase-locked spatio-temporal (and velocity) harmonics are created, even though only one spatio-temporal harmonic was driving the plasma and only for a finite time duration, turned on and off adiabatically. These self-consistent field structures that arise can trap, untrap and retrap particles whose orbits are near the non-stationary large-excursion separatrix. These lost particles are not retrapped where they were released, but elsewhere in distant troughs of the field, after many oscillations in either direction. While KEEN waves are harder to excite and to self-organize, they are very robust to perturbations. On the other hand, delicate resonant modes have the opposite property. Namely, they are readily excited (even via infinitesimal amplitude perturbations) but at larger amplitudes, they are easily detuned and deformed.  

The self-adjusting multi-harmonic field structures in KEEN waves can trap enough of the particles to maintain themselves in perpetuity within the Vlasov-Poisson or Vlasov-Maxwell \cite{bedros_ref1} system of coupled integro-differential equations. With KEEN waves, deeply trapped particles remain trapped, but the separatrix regions harbor far more complicated dynamics. The weaker the drive, the more the entire fate of the mode is dictated by proliferating separatrix regions around small vorticlets. The stronger the drive, the deeper the wells and the less stringent this loss mechanism is on the overall sustainability of the mode. But as the drive amplitude is diminished, the dynamics becomes fragmentary and small vortices may or may not be able to remerge and form a KEEN wave. The process takes longer and longer at smaller and smaller drives and may become disrupted in reality by other physical processes such as collisions or side-losses, which are not included in this $1D\times1D$ Vlasov-Poisson model. 

We have observed KEEN waves in the laboratory driven exactly as stated above, via the ponderomotive force emanating from the optical mixing of two laser beams. For details of those experiments conducted on the Trident laser system at Los Alamos, see \cite{kline}. In such experiments, $3D$ effects, ion dynamics and all manner of collisions were of course inevitable. This gives us assurance that the existence of KEEN waves is not strictly limited to the low dimensional numerical investigations we have led so far. Similarly, Vlasov-Maxwell simulations have shown that KEEN waves can be driven and sustained with interesting interactions with electron plasma waves, see \cite{bedros_ref1}

Here we wish to demonstrate how to perform efficient simulates of  the transition region between a well formed KEEN wave and more fledgling scenarios, by varying the drive time of the ponderomotive force, and following the long time evolution of the structures that remain and persist, long after the drive is turned off (and the system has a time independent Hamiltonian). Two crucial new elements are introduced to make these advances possible. First a new Vlasov solver is implemented which does not use a fixed velocity grid, while still using a fixed spatial grid. Instead, a two velocity-grid system is adopted. Fine scale resolution is maintained around the phase velocity of the driven wave, where the vortical structures are formed. This is surrounded by an external velocity region where only coarser scales are resolved. The external region reacts to KEEN waves by absorbing its wake and this must be computed correctly in order to gather the correct overall charge density to be used as the source term of the Poisson equation. The variable grid technique allows the inner $v$-space region to shrink in absolute size while maintaining high resolution within, without requiring massive fine scale resolution in the exterior velocity regions. Just the perturbed region is finely resolved, no matter how small it gets. Secondly, since long time evolution of slowly forming objects in phase space are sought, it is highly desirable to have a way to speed up the calculations without losing accuracy, especially when faced with proliferating small phase space structures with many separatrix regions, which may be eventually merging, that have to be tracked. This is the dilema in low amplitude drive cases unless a variable velocity grid and high order time splitting methods comes to the rescue. 
 
The mutual attraction and commingling of many vorticlets at low amplitude drive can cause one to question the accuracy and fidelity of low-order, long time simulations. The usual modern answer, when facing such nonlinear dynamical problems, is to appeal to symplectic and higher order integration schemes. Simplifications in the higher order commutator evaluations are found that ease the computational burden \cite{ccfm20pp}. In fact, we adopt an easy to implement sixth order operator-splitting scheme, with minimal additional computational burden. This in turn allows large time steps to be taken (such as 0.5 or 0.25) without appreciable loss of accuracy over lower order methods with much finer time steps. 

Armed with these two new techniques, we have simulated weakly driven KEEN waves with a sequence of drive durations that extend from no apparent KEEN wave formation all the way to the canonical KEEN wave conditions reported previously \cite{bedros04,vog,bedros14b}. 

Intuitively, we know that it is the amount of energy that is directly coupled into the KEEN wave that will dictate its final size in phase space or its peak electric field strength. But the scaling with drive time is not linear. There is, in fact, a saturation that sets in of the directly driven mode  (or the first harmonic). The driven mode at its maximum amplitude will collapse and channel its energy to higher order modes, starting with the second harmonic. This is repeated between the second and the third harmonics and so on, depending on the strength of the drive and the responsiveness of the plasma (more harmonics are implicated when driven at longer wavelengths, than at shorter wavelengths, generally. Many harmonics correspond to small values of the wavenumber of the drive, measured in inverse Debye lengths). We excite KEEN waves for a finite amount of time after which the drive is shut off and the modes rearrange themselves and adjust to a trapped-untrapped and retrapped distribution of particles that restore energy to the field, maintaining a coherent non-decaying field structure. An important question is how long must the drive last, so that trapping sets in for large scale vortices, that do not disintegrate after the drive is turned off, but which can maintain the KEEN wave?

A simplified physical picture of KEEN wave formation is as follows. The wavenumber of the sinusoidal drive sets a length scale for self-organization via trapping. The amplitude and duration of the drive, on the other hand, set a different rule relating the particle energy to the trapping period. The typical scaling of range of velocities affected by a driver field of fixed amplitude, duration and wavenumber is well known. Assuming a harmonic field structure, the scaling of vortex width in velocity is as the square root of the response field, which itself is proportional to the maximum of the driver electric field amplitude. The proportionality constant is here denoted by $\alpha$. What we need is a self-sustaining structure in phase space which is rotating and also translating at the phase velocity of the drive. Rotation or trapping oscillations imply a fixed energy orbit, $E$. That scaling is, in normalized units explained below:

$$  |v_{max}| = \sqrt{2 \big[ E + \alpha \, a_{\rm Dr} \big] } \,\, . $$

In this simple picture, we see that if the drive's amplitude or duration is not large enough, the velocity range over which the fledgling vortex will be forming may be  far smaller than the wavelength of the drive (in natural units for Vlasov-Poisson system, as explained below). Then, if the drive is turned off, a single well of the size of the driver wavelength will not be able to trap a sufficient number of particles to sustain itself. There are two inherent scales which are not commensurate, in that case. This will led to immediate break up of the (seemingly) long wavelength structure to a set of shorter wavelength ones which are of the same length scale as the velocity disturbance achieved. This is the collapse of the input energy into small vorticlets of the size of the square root of the response (self-consistent) electric field. Local trapping can occur but on a much shorter length scale than the wavelength of the applied field. However, if the drive amplitude and its duration are long enough to make the maximum bounce velocity of a large number of particles, in thermal speed units, be the same size or larger than the wavelength of the drive in Debye length units, then we will be able to sustain self-organized KEEN waves in perpetuity (within the confines of the equations making up the Vlasov-Poisson or Vlasov-Maxwell models). Higher harmonics will lead to the periodic shedding of smaller vortices which will remerge eventually and contribute to the non-stationarity. 

In these more strongly driven cases, the KEEN wave formed will involve multiple higher harmonics which allow the overall trapping of shallower trapped particles by sometimes shedding small vorticlets which eventually merge in with the main one. Stretching and folding around the separatrices is the active mechanism leading to chaotic orbits. When the drive amplitude is large, a large vortex dominates and all smaller vortices formed can be swallowed up by the large one. When the drive is weak, a large number of small vorticlets form at the end of the drive period and they may or may not be able to coalesce into a deep enough trapped particle potential at or near the scale of  the original driving wavelength. These two scenarios have their own chaotic properties that are dissimilar. They put demands on phase space resolution in different ways. This will be discussed in this paper where cases in both regimes of drive are shown in detail. In the weak drive case, after the drive is turned off, one can only rely on a merging of the small vorticlets for a KEEN wave to eventually form at or near the original wavelength scale. This is neither an efficient nor a reliable process. The weaker and shorter lived the drive, the less chance there is of this process succeeding quickly enough. Thus a (complicated)  threshold does exist for true KEEN wave self-organization to occur.  Very reliable simulations are required in order to capture this physics with confidence. 

In this paper, we show a number of examples of KEEN wave generation and evolution and compare uniform velocity grid simulations to two grid, variable resolution simulations, with the added feature of an easy to implement sixth-order, symplectic, time-splitting integrator. These features will allow us to observe the birth pains of KEEN waves as they transition between vorticlets attempting to merge, and full formed, self-organized multiple phase-locked harmonic structures shedding small vortices periodically and swooping them back up in chaotic fashion. We will show a series of cases at different drive amplitudes and for different drive durations.

The $1D\times1D$ Vlasov-Poisson equations for the electron distribution function $f=f(x,v,t)$ with an external  ponderomotive force drive electric field $E_{\rm Pond}$ reads \cite{bedros04,vog,bedros14b}:
\begin{eqnarray*}
\partial_t f + v\partial_x f+(E-E_{\rm Pond})\partial_v f = 0,\ \partial_x E = \int_{\mathbb{R}} f dv -1. 
\end{eqnarray*}
Here, space is measured in Debye length units, velocity in thermal speed units and time in inverse electron plasma frequency units. $E_{\rm Pond}(x,t)$ may be the sum of elements each of which is of the form
$$
E_{\rm Pond}(x,t) = a_{\rm Dr}\,k_{\rm Dr}\,a(t) \,\sin(k_{\rm Dr}\,x-\omega_{\rm Dr} \,t), 
$$
with the adiabatic switch on and switch off form being given by $$
a(t)=\frac{g(t)-g(t_0)}{1-g(t_0)},
$$
$$
 g(t) = 0.5(\tanh(\frac{t-t_L}{t_{wL}}) - \tanh(\frac{t-t_R}{t_{wR}})),
$$ 
$t_0=0,\ t_L=69,\ t_{wL}=t_{wR}=20,\ t_R=207+T_{\rm Dr}$, 
 $k_{\rm Dr} = 0.26$, $\omega_{\rm Dr} = 0.37$.
The initial condition is a spatially uniform Maxwellian velocity distribution at a fixed temperature. 
$$
f_0(x, v) = \frac{1}{\sqrt{2\pi}}\exp\left(-\frac{v^2}{2}\right), \;\;\; (x,v) \in [0,2\pi/k_{\rm Dr}] \times [-6,6].   
$$
Two very distinct cases are  considered:{\bf Canonical drive}, with $T_{Dr}=100$ and $a_{\rm Dr}= 0.2$, and
{\bf Weak drive}, with $T_{\rm Dr}=200$ and $a_{\rm Dr}=0.00625$. The initial results we will show are for relatively short total run time simulations where a dozen or so trapping times have elapsed. This corresponds to a time, in plasma frequency inverse units, of $T=1000$ in the canonical case
and $T=5000$ in the weak drive case. Then, by varying the drive duration $T_{\rm Dr}$, in the weakly driven case, $a_{\rm Dr}=0.00625$,
we will capture the complex fragmentation and merger dynamics of vorticlets leading to the formation of KEEN waves, whenever that occurs. Chasing KEEN wave formation under such circumstances will require simulations that are at least 10 times longer. This will bring into question fine scale fidelity which can not be guaranteed. We trust that the global features will be correct however. 


\noindent Simulating KEEN wave can be challenging, as we have explained above.  Resolving small localized phase space structures consistently, and doing so to large times, typically requires a large number of phase space grid points and the use of small time steps. Achieving high order accuracy in both (space and velocity direction) interpolation and in time-splitting strategies have been discussed recently \cite{cgm2012,vog,steiner_vlasovia2013}.
\noindent Here, we continue the work began in \cite{vog}, where the canonical KEEN wave case was simulated with a number of arbitrary order interpolation techniques than in the cubic spline code that led to the discovery and initial exhaustive studies of KEEN wave dynamics first reported in \cite{bedros04}. These semi-Lagrangian methods were equipped with high order  spline, Lagrange or Hermite interpolation, and used Strang splitting. So while the interpolation order was being varied at will, the temporal resolution was never keeping pace leading to surprising results such as more uncontrolled fluctuations 
The effects of single vs double precision were studied and a $\delta f$ method was adopted in order to improve the reach of single precision simulations on GPUs. Without that $\delta f$ approach, numerical noise would contaminate the results too easily. The use of GPU permitted run efficiency in finely resolved simulations ( we used $2048\times2048$ phase space grids in double precision and $4096\times 4096$ in single precision; with a very small time step, such as $0.01$). However, we were not able to achieve late time convergence of the full dynamics of even the first few spatial Fourier modes of the charge density. While the phase space vortical structures making up KEEN waves looked quite similar, there were features with detectable deviations out to late times, which changed at different resolutions.  
The simulations we will perform here require high resolution around the phase velocity of the driven wave, and this region becomes smaller and smaller when the drive amplitude is diminished.
We have performed many runs at different choices of numerical parameters. We report here on the results of parameter choices which worked well.
This should facilitate future choices of suitable numerical parameters. These constitute but a first guiding step for a more rigorous and physically motivated parametric study \cite{bedros14b}. 
More elaborate physically motivated diagnostics will also be deployed in those KEEN wave dynamics studies. As for this paper, in Section \ref{sec:1}, we
give a detailed account of the main ingredients of the numerical method.  At the center of it is the conservative, non-uniform, cubic splines reconstruction method and the high order splitting scheme in time with cancellation tricks that lead to fast computation.  
Section \ref{sec:2} is devoted to the elaboration of numerical results: first, two distinguishing limiting drive amplitude cases will be discussed in detail. One is the canonical case, while the other is the weak drive case.  We will also vary the drive time, for the weak drive amplitude case, to better capture the complex KEEN wave creation dynamics. This will allow us to observe vorticlet formation and their possible merger dynamics. 




\section{New Numerical Methods}
\label{sec:1}
Obtaining point-wise numerical convergence is not necessary nor expected in highly chaotic, large dimensional dynamical systems. All we can hope for are consistent average quantities being tracked as resolution in space, velocity and time grids are increased. Here, we add three new ingredients beyond the high order scheme advances reported in \cite{vog}. The first is the MPI parallelization of the code (previous versions were OPENMP or using GPU). The new code is developed within the framework of the (Semi-Lagrangian) Selalib library \cite{selalib}. The second, is the use of non-uniform, conservative, cubic-splines, in order to capture the intricate internal dynamics of KEEN wave vorticlet merger arising in smaller and smaller velocity ranges, as drive amplitude is diminished. The third, is the use of sixth order time splitting \cite{ccfm20pp}, in order to be able to use larger time steps, maintain high temporal accuracy and be confident while integrating to large times. We describe the two numerical method advances below, but do not discuss parallelization strategies further. 

\subsection{The variable resolution, two-grid mesh in velocity}

\noindent We generate a $1D$ velocity mesh that is finely resolved in a (relatively) small region, while maintaining a coarser grid elsewhere. This is because for typical drive amplitudes and durations, the distribution function is perturbed only in the immediate vicinity of the phase velocity of the driving field  
 and scales roughly as the square root of the drive amplitude. The strategy adopted here defines a coarser grid everywhere which has a refined interior region.
This grid construction has the advantage of simplicity. The mesh spacing used on the coarse and fine grids are: 
$$\Delta v_{\rm coarse}=\frac{v_{\rm max}-v_{\rm min}}{N_{\rm coarse}},\ \Delta v_{\rm fine}=\frac{v_{\rm max}-v_{\rm min}}{N_{\rm fine}}$$ 
and $N_{\rm fine}$ is an integer multiple of $N_{\rm coarse}$.

\noindent The refined zone is chosen with $0\le i_1<i_2\le N_{\rm coarse}$
and the total number of cells is
$$
N=i_1+N_f+N_{\rm coarse}-i_2,\ N_f=\frac{N_{\rm fine}}{N_{\rm coarse}}(i_2-i_1)
$$

\begin{center}
\begin{tikzpicture}[scale=1]
\foreach \i in {0,...,5}
{
\path (0.4*\i,0) coordinate (X\i);
\fill (X\i) circle (1pt);
}
\textcolor{red}{
\foreach \i in {1,...,20}
{
\path (2+0.1*\i,0) coordinate (X\i);
\fill (X\i) circle (1pt);
}
}
\foreach \i in {1,...,5}
{
\path (4+0.4*\i,0) coordinate (X\i);
\fill (X\i) circle (1pt);
}

\textcolor{blue}{
\foreach \i in {0,...,15}
{
\path (0.4*\i,0) coordinate (X\i);
\fill (X\i) circle (1pt);
}
}

\textcolor{violet}{
\foreach \i in {15,...,15}
{
\path (0.4*\i,0) coordinate (X\i);
\fill (X\i) circle (2pt);
}
}

\textcolor{violet}{
\foreach \i in {0,...,0}
{
\path (0.4*\i,0) coordinate (X\i);
\fill (X\i) circle (2pt);
}
}

%
\path (3.,0.6) node (v1) {\textcolor{red}{$N_f$}};
\path (1.,0.6) node (v1) {\textcolor{blue}{$i_1$}};
\path (5.,0.6) node (v1) {\textcolor{blue}{$N_{\rm coarse}-i_2$}};
\path (0.5,-0.3) node (xi) {\textcolor{violet}{$v_0=v_{\rm min}$}}; 
\path (2,-0.3) node (xi) {$v_{i_1}$}; 
\path (4,-0.3) node (xi) {$v_{i_1+N_f}$}; 
\path (6.8,-0.3) node (xi) {\textcolor{violet}{$v_N=v_{\rm max}$}}; 
\path (0,0) coordinate (vmin);
\path (6,0) coordinate (vmax);
\draw (vmin) -- (vmax);
%
%
%
%
%
%
\end{tikzpicture}
\end{center}
In order to generate such a mesh, we have chosen as input parameters the total number of cells $N$, the region where we want to refine
$v_{\rm min} \le a < b \le v_{\rm max}$ and an integer ratio $r \not=1 $. From these values, we look for 
$$i_1,\ i_2,\ N_{\rm coarse} $$
such that
$$
v_{\rm i_1} \simeq a,\ v_{\rm i_1+N_f} \simeq b,\ r = \frac{N_{\rm fine}}{N_{\rm coarse}}.
$$
The algorithm we have adopted is this. We first write 
$$
\alpha = \frac{a-v_{\rm min}}{v_{\rm max}-v_{\rm min}},\ \beta = \frac{b-v_{\rm min}}{v_{\rm max}-v_{\rm min}}
$$
and compute
$$
N^*_{\rm coarse} = \left\lfloor \frac{N}{1+(\beta-\alpha)(r-1)}\right\rfloor,\ N_{\rm fine} = rN^*_{\rm coarse},
$$
together with
$$
i_1^* = \left\lfloor \alpha N^*_{\rm coarse}\right\rfloor,\ \ell =  \left\lfloor\frac{N-N^*_{\rm coarse}}{r-1} \right\rfloor,
$$
in order to obtain
$$
N_{\rm coarse} = N-\ell (r-1),\ i_1 =  \left\lfloor \alpha N_{\rm coarse}\right\rfloor,\ i_2 = i_1+\ell.
$$
Other strategies to be pursued may be more efficient, including changing the grid smoothly in order to absorb the shock of a suddenly changed mesh spacing, which here goes directly from
$\Delta v_{\rm coarse}$ to $\Delta v_{\rm fine}$.

\subsection{Conservative cubic splines on a non-uniform mesh}

\noindent One popular semi-Lagrangian method for the numerical solution of the Vlasov-Poisson equations is dimensional splitting \cite{cheng} while using cubic splines for the interpolation. In this paper, the mesh in velocity is not uniform so that the adoption of non-uniform cubic splines is required. However, we lose the conservation of mass property of the algorithm, when the mesh is non-uniform. That can have deleterious effects on the numerical results of KEEN wave formation, as we have observed. For instance, in \cite{vog}, we observed the devastating effects of non-conservation of mass due to the use of single precision computations.
Thus, instead of looking at the classical advective form of the constant advection equation, we can consider the conservative form \cite{cms}. That means that we reconstruct the primitive function, using the same interpolation operator. We have to be careful about the boundary conditions, which here we considered to be periodic. We will have also to shift the unknowns to the middle of the velocity cells, in this non uniform setting (this is not necessary, for the uniform grid case and edge definitions work). We now show in full detail the different steps of the new algorithm.

\noindent Thanks to dimensional splitting, we are lead to solve
$$
\partial_t u+c\partial_vu = 0,
$$
over a  time step $\Delta t$ (or a fraction of a time step, in fact), with unknowns 
$$
u_{j+1/2}(t) = \frac{1}{v_{j+1}-v_{j}}\int_{v_{j}}^{v_{j+1}}u(v,t)dv,\ j=0,\dots,N-1.
$$ 
That is, we are supposed to know $$u_{j+1/2}^{\rm old} \simeq u_{j+1/2}(0),\ j=0,\dots,N-1,$$ and we want to compute
$$
u_{j+1/2}^{\rm new} \simeq u_{j+1/2}(\Delta t),\ j=0,\dots,N-1.
$$
\noindent Using the conservation of volume, we have the relation
$$
\int_{v_{j}}^{v_{j+1}}u(v,\Delta t)dv = \int_{v_{j}-c\Delta t}^{v_{j+1}-c\Delta t}u(v,0)dv.
$$
\noindent We first compute 
$$
U_j = \sum_{k=0}^{j-1}(u^{\rm old}_{k+1/2}-M)(v_{k+1}-v_k),\ j=0,\dots,N,
$$
with
$$
M = \sum_{k=0}^{N-1}u^{\rm old}_{k+1/2}(v_{k+1}-v_{k}).
$$
Note that, by construction, we have $U_N=U_0=0$.
We then define the non-uniform cubic spline interpolation of the primitive, that is the unique piecewise cubic polynomial function $U_h\in C^2_{\rm per}(v_0,v_N)$
 satisfying
$$
U_h(v_j) = U_j,\ U_h \textrm{ polynomial on } \ [v_j,v_{j+1}],\ j=0,\dots,N-1. 
$$
This can be done classically by solving a system which is almost tridiagonal, for computing the spline coefficients or the Hermite derivatives; see e.g. \cite{cygne}.
\noindent Note that the primitive is periodic, thanks to the choice of the integration constant $M$.
\noindent Finally, we compute
$$
U_j^{\rm new} = U_h(v_j-c\Delta t),\ j=0,\dots,N-1,
$$
and get the unknowns updated by
$$
u^{\rm new}_{j+1/2}= \frac{U_{j+1}^{\rm new}-U_{j}^{\rm new}}{v_{j+1}-v_j}+M,\ j=0,\dots,N-1.
$$
Note that the method is conservative by construction, as we get
$$
\sum_{j=0}^{N-1}u^{\rm new}_{j+1/2}(v_{j+1}-v_{j}) = M.
$$
Other strategies can be envisaged, by using the special structure of the two-grid mesh \cite{posterkeen}, but are not further developed here. 


\subsection{High order time splitting}

We fix the time step $\Delta t$ and consider a list of coefficients $a_1,\dots,a_s$, with $s\in \mathbb{N}^*$, together with a coefficient $\sigma_{\rm init}\in\{0,1\}$.


\noindent For $n\in \mathbb{N}$, we know 
$$
f^n_{i,j}\simeq f(n\Delta t,x_i,v_{j+1/2})
$$
\noindent Index $i$ will go from $0$ to $N_x-1$ and $j$ from $0$ to $N_v-1$.
\noindent We fix $t^*\rightarrow n\Delta t$ and $\sigma\rightarrow \sigma_{\rm init}$, and start with $f^*_{i,j}\rightarrow f^n_{i,j}$. 

\noindent For each $k=1,\dots,s$, we perform the $\mathcal{T}^\sigma$ advection  over a time step $\Delta \tau \rightarrow a_k\Delta t$ and then update
$\sigma\rightarrow 1-\sigma$.

\noindent Here $\mathcal{T}^0$ (advection in $x$) consists in solving over a substep $\Delta \tau$
$$
\partial_tf(t,\cdot,v_{j+1/2})+v_j\partial_xf(t,\cdot,v_{j+1/2}) = 0,\ 
$$
to update $f^*_{i,j}$. At the end, we update $t^*\rightarrow t^*+a_k\Delta t$.

\noindent $\mathcal{T}^1$ (advection in $v$) consists in computing the electric field $E^*(x_i)$ via the Poisson equation (see e.g. \cite{vog}) and solve over a substep $\Delta \tau$

$$
\partial_tf(t,x_i,\cdot)+(E^*(x_i)-E_{\rm app}(t^*,x_i))\partial_vf(t,x_i,\cdot) = 0,
$$ 
to update $f^*_{i,j}$.
At the end of the substep $s$, we get $f^{n+1}_{i,j}\rightarrow f^*_{i,j}$.

\noindent Classical Strang splitting, that will be used here for comparison, corresponds to $s=3,\ a_1=1/2,\ a_2=1,\ a_3=1/2$ and $\sigma_{\rm init}=1$.

\noindent We have developed new efficient high order schemes for Vlasov-Poisson, see \cite{ccfm20pp},
exploiting the specific 
structure of the Vlasov-Poisson system of equations in $1 D$. The $6$th order Vlasov-Poisson splitting scheme, that is used to produce the numerical results to be shown below,
has the following coefficients: $s=11$,
$$
\begin{array}{ll}
a_1= &0.0490864609761162454914412\\
&-2\Delta t^2(0.0000697287150553050840999),\\
a_2 = &0.1687359505634374224481957,\\
a_3 = &0.2641776098889767002001462\\
&-2\Delta t^2(0.000625704827430047189169)\\
&+4\Delta t^4(-2.91660045768984781644\cdot 10^{-6}),\\
a_4= &0.377851589220928303880766,\\
a_5 = &0.1867359291349070543084126\\
&-2\Delta t^2(0.00221308512404532556163)\\
&+4\Delta t^4(0.0000304848026170003878868)\\
&-8\Delta t^6(4.98554938787506812159\cdot 10^{-7}),\\
a_6 = & -0.0931750795687314526579244,\\
\end{array}
$$
together with $a_{6+i}=a_{6-i},\ i=1,\dots,5$ and $\sigma_{\rm init}=1$.

\section{Numerical results}
\label{sec:2}

\noindent In all the numerical results shown below, we have used $v_{\rm max} = 6$, and a Lagrange interpolation of degree $17$ in $x$ (see \cite{vog}). In the $v$-direction, 
we have used conservative, non-uniform cubic splines. When the mesh is uniform, this corresponds to classical cubic splines. For the sake of simplicity, we have used the same code
for the uniform and non-uniform meshes. This can affect speed comparisons, as we would normally expect faster code performance for uniform cubic splines. We should stress that much code optimization remains to be done in the future in both the uniform and the non-uniform mesh cases. 
For non-uniform meshes, we use an integer ratio $r=32$, $a=1.2,\ b=1.6$, for the weakly driven, and $a=0.375,\ Êb=2.25$ for the canonical drive case.
Simulations are run on Hydra (the computing centre of the Max Planck Society)  and Helios Computational Simulation Centre, International Fusion Energy Research Centre of the ITER Broader Approach) supercomputers, with typical runs on $256$ processors ($16$ nodes; each node having $16$ threads); more specifically, on helios, we use the selavlas allocation on the main partition or the mic\_eu allocation, which
is a partition dedicated to accelerated nodes, which contains also standard nodes, that we used. 
The parameters of the different runs reported on here are summarized in Table \ref{table1}.

\subsection{Efficiency of the code}
In order to measure the efficiency of the code, which is useful for future optimization and for comparisons with earlier results, we define the efficiency by
$$\textrm{eff}=(s\cdot N_x\cdot N_v\cdot T/\Delta t)/(\textrm{time}\cdot10^{6}\cdot \textrm{proc}),$$
where $s$ is the number of substeps per iteration ($s=3$ for Strang splitting and $s=11$ for the $6$-th order splitting), proc is the number of processors, time is the duration of the simulations in seconds.
This measure is a good indicator of the efficiency of implementation of the scheme. In order to compare with other results (such as those in \cite{vog}), we can obtain an efficiency of 
$29$ (resp. $15$) on a $256\times 256$ (resp. $2048\times 2048$) grid on a sequential code. Using GPU and optimized fft, we obtained an efficiency of 
$73$ (resp. $414$) on a $256\times 256$ (resp. $2048\times 2048$) grid. So, we can remark that GPU is not more competitive, once the number of processors exceeds $128$, that is to say, past $8$ nodes.
There is still ample room for improvements in the implementation. One challenge is to use multi-GPU or MIC nodes; another, is to improve the efficiency of the MPI parallelized code. One step along these lines would be to use a mixed OPENMP/MPI parallelization.

\subsection{Diagnostics}
\noindent We resort to a very limited set of diagnostics in this paper which show how well the code performs and to give but a glimpse at the physics of KEEN waves and their formation. The diagnostics we show are:  
\begin{itemize}
\item $\delta f$ $= f-f_{0},$ the change in the distribution function with respect to the initial condition, as a function of time. 
\item the $5$ first $\rho$ harmonics. That is the absolute values of the
Fourier mode amplitudes of $\rho = \int_\mathbb{R} f dv$, from mode $k = 1$ to mode $k = 5$
\item The relative $L^2$-norm error (exact $L^2$-nom is $(\frac{\sqrt{\pi}}{0.26})^{1/2}$).
\item The time evolution of  the $L^2$-norm in velocity of spatial Fourier modes of $\delta f:$  $$\sqrt{\int_\mathbb{R}|\hat{f}^{(k)}(v,t)|^2dv}$$ for modes $k=1,2,\dots,6$ and $8,\ 12,\ 16,\ 20$. This diagnostic shows KEEN wave formation to correspond to specific set separations between the amplitudes of the low order harmonics (say modes 2, 4) being established and maintained indefinitely. When all mode amplitudes are small, fluctuating and indistinguishable in relative amplitude then the resulting phase space structure is not yet coalescing or self-organizing into a clean KEEN wave. This is a velocity space RMS diagnostic of each spatial Fourier mode vs time. 
\end{itemize}

%
%
%
%

\subsection{The weakly driven case}
\begin{figure}[h!]
\begin{tabular}{c}
\includegraphics[trim=4.cm 4cm 2cm 2.4cm, clip=true, width=0.45\linewidth,angle=0,height=7.5cm]{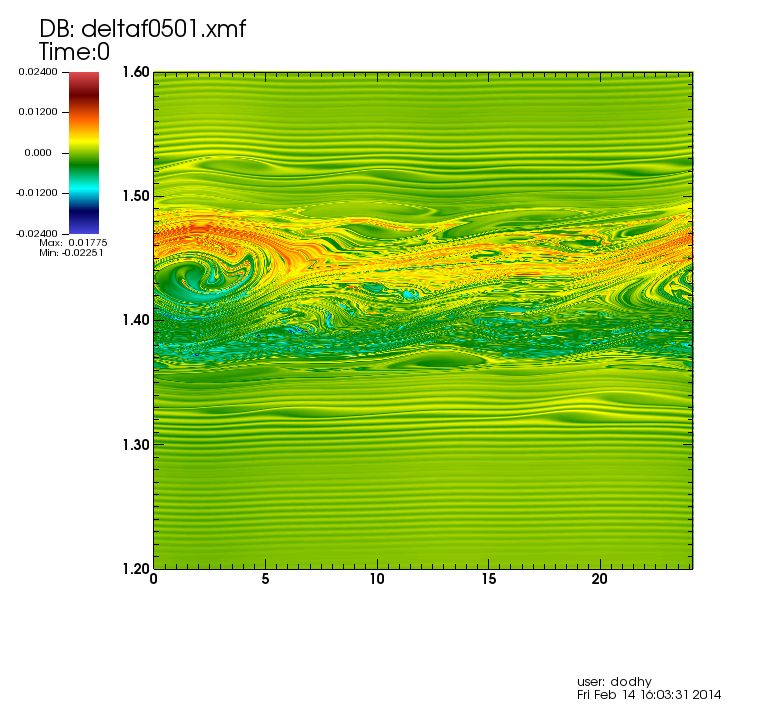}
\\
\includegraphics[width=0.5\linewidth,height=4cm]{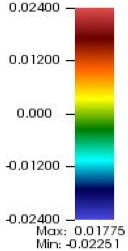}\\
\includegraphics[trim=1.cm 0.5cm 0cm 2cm, clip=true, width=0.45\linewidth,angle=0,height=5.5cm]{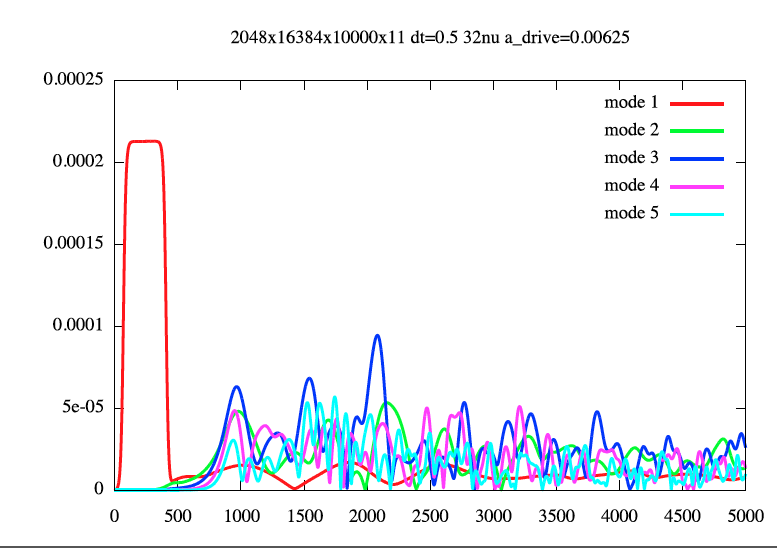} 
\end{tabular}
\caption{$\delta f$ distribution function $f-f_0$ at time $T=5000$ as a function of $(x,v) \in [0,4\pi]\times [1.2,1.6]$ (top) and time evolution of the amplitudes of the first $5$ $\rho$ harmonics (bottom), 
in the {\bf small drive} amplitude case: canonical run with $6$th order time scheme and non-uniform velocity mesh.
Parameters are $N_x=2048,\ N_v=16384,\ \Delta t=0.5$ non-uniform, $6$-th order time scheme (run13).}
\label{fig1_small}
\end{figure}

We consider at first the small drive amplitude case, that is with drive amplitude $a_{\rm Dr}=0.00625$ and drive time $T_{\rm Dr}=200$.
Results are given on Figures \ref{fig1_small} and  \ref{Tdrive200}. 

\noindent In the beginning, no KEEN wave is yet forming and filamentation occurs
(as for the non linear Landau damping test problem). At time $T=5000$, we can detect a large number of vorticlets with one larger one swallowing up its neighbors; looking further in time
(see Figure  \ref{Tdrive200}), we see that the KEEN wave persists.

\noindent Filamentation is mainly along the $v$-direction and thus high velocity resolution is needed. Note that $N_v=16,384$ in the inner region would correspond to $N_{\rm fine}=258,432$ points in velocity, had the mesh been uniform in all of v. In our case,  only the region $[1.2,1.6]\subset [-6,6]$ is refined, which represents $1/30$ of the whole velocity region. 
A non-uniform grid code is imminent reasonable in this case. 
In comparison to the canonical drive case, we have to wait much longer until vortical merging structures develop sufficiently. High order time splitting permits error control in time for a fixed, not too small a time step. Here we take $\Delta t=0.5$; which permits us to perform a simulation with not too many time steps (here $20000$ for reaching the time $T=10000$).

\subsection{The canonical drive case}
\begin{figure}[h!]
\begin{tabular}{c}
\includegraphics[trim=4.5cm 4cm 2cm 2.4cm, clip=true, width=\linewidth,angle=0,height=7.5cm]{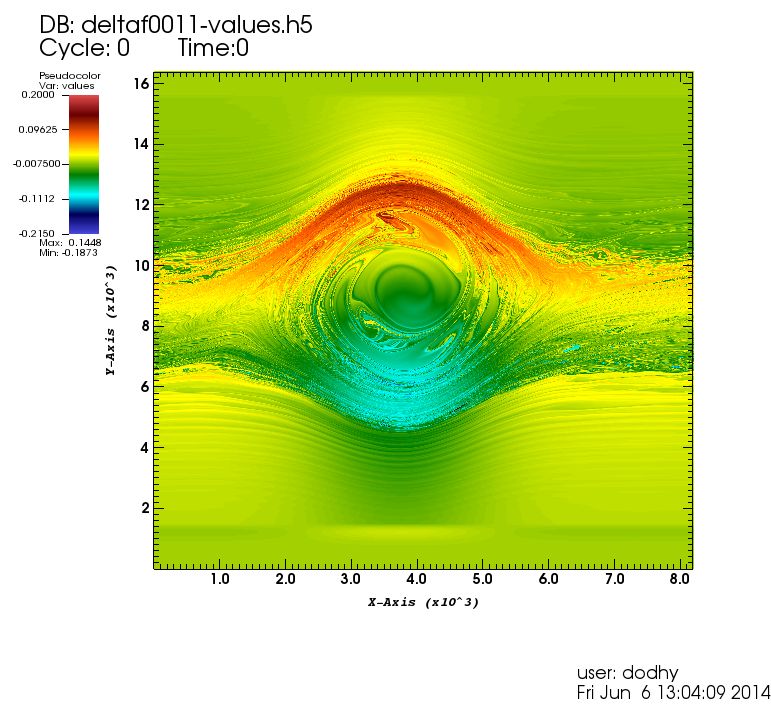}\\
\includegraphics[trim=0.cm 0cm 0cm 0.2cm, clip=true,width=0.5\linewidth,height=4cm]{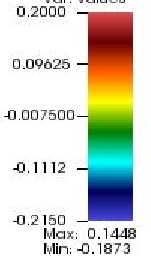}\\
\includegraphics[trim=1.cm 0.5cm 0cm 2cm, clip=true, width=0.5\linewidth,angle=0,height=5.5cm]{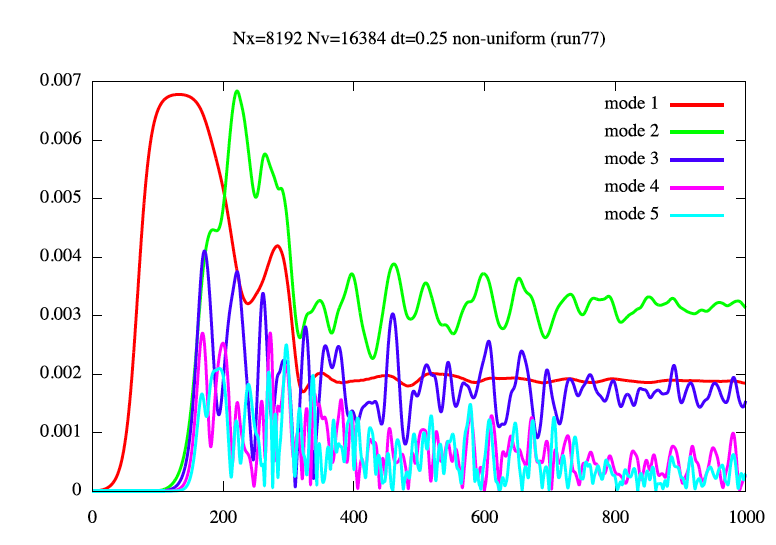} 
\end{tabular}
\caption{$\delta f$ distribution function $(f-f_0)(x_i,v_j)$ at time $T=1000$ as a function of $(i,j) \in [0,N_x]\times [0,N_v]$ (top) and time evolution of the first $5$ $\rho$ harmonics (bottom), 
in the {\bf canonical drive} case: canonical run with $6$th order time scheme and non-uniform velocity mesh.
Parameters are $N_x=8192,\ N_v=16384,\ \Delta t =0.25$ non-uniform, $6$-th order time scheme (run77).}
\label{fig2_std}
\end{figure}

We return to the canonical drive case, that is with drive amplitude $a_{\rm Dr}=0.2$ and drive time $T_{\rm Dr}=100$.
Results are given in Figure \ref{fig2_std}. For this case, we already know what to expect with sufficient confidence \cite{bedros04,bedros14b,vog}. The benefit of the non-uniform mesh
is lower here, but still present. The refined region is now $[0.375,2.25]\subset [-6,6]$, that is $1/6.4$ of the whole velocity region. 
KEEN waves develop earlier and vortex mergers occur far more rapidly and decisively. Without high enough spatial and velocity resolution, excessive numerical diffusion occurs which smooths the distribution function unduly and prematurely. The higher the resolution, the later the noticeable effects of the numerical diffusion process appear. In order to have "converged" results, we use $N_x=8192$ and $N_v=16384$.
An equivalent uniform velocity mesh would require $N_{\rm fine}=89792$ points. With so many points, we would even have difficulty visualizing the entire distribution function. Instead, we plot the distribution function just as
$f(x_i,v_j),$ a function of $i$ and $j$. This permits ease of visualization, by removing the full mesh information. Nevertheless, we can still see the transition between the coarse and fine meshes and the extent of the fine mesh in comparison to the coarse one. Note that in zoomed in plots, the coarse region is not plotted; we always zoom in on the KEEN wave region.

\subsection{First remarks regarding convergence}
We focus our attention on the canonical drive case, where convergence was not achieved in previous studies using low or high order interpolation schemes. To move beyond that, we have relied upon the convergence rates of a number of global diagnostics to guide us  in refining velocity, space and/or time resolution. As we have explained, small vortices shed from the main large vortex street of a KEEN wave, if not tracked properly, will cause fluctuations in the density response in time which will change with changing resolution. With enough spatial and velocity resolution this problem will be diminished. Since at healthy drive levels, there will be small vortex shedding and remerger on many small scales (with respect to the drive wavelength), the need for very high resolution in phase space is expected, if fine details are to be tracked for long times.  

\noindent Since fine scale self-organization is being created via the large (drive wavelength) scale KEEN wave in the canonical case. While in the weakly driven case, the initial formation process has not survived the drive and so smaller scales than the initial fragmented scale vortices are not being created for lack of larger scale self-organization to drive them. For the canonical drive case, we had to use a time step of $\Delta t=0.25$ even for a sixth order scheme in order to track the small phase space vortical structures correctly.
The use of the $6$th order in time scheme really kicked in and showed its power only after enough resolution was being achieved in the active region in phase space in space and velocity. When dealing with lower resolution runs, we found our results to be much more sensitive to the time step used. There is an important lesson here which pertains to method of Vlasov equation solution used. When interpolation schemes are deployed, such as splines, 
too many repeated applications, as necessitated by very small time steps, can corrupt the solution with little real benefit. Therefore, seemingly paradoxically,  the best result for a fixed phase space resolution will not be achieved by taking the smallest possible time step (in the context of semi-Lagrangian schemes). Results  in $\rho$-harmonics and the $L_2$-norm diagnostics amply confirm this observation. It is difficult to call out the proper time step to use a priori. For KEEN waves, the difficulty can be traced back to the fact that many harmonics (spatial scales) are involved in the KEEN wave creation process, followed by a few important scales once the KEEN wave has been formed. And yet, that overall number is still in the tens of modes compared to the driver wavelength. During the creation process, hundreds of harmonics share energy and interact. The retrenchment of coherent energy back down to a few modes is a telltale sign of self-organization. For weak drive, this may proceed differently. There, immediate fragmentation of vortical structures occurs once the drive is shut off and the regrouping or remerged between these small vortices may or may not occur in finite time. While the inherent time scale of evolution is slower, since there is less energy injected into the system to cause violent behavior, the phase space complexity is more. So there are two completely opposite scenarios. In weak drive, slowly evolving complex structures must be tracked in phase space. While for large enough drive, very fast and energetic shredding and remerged of small pieces onto a massive core that comprises a KEEN wave must be monitored to make sure they do not cause structural instability. Also, when two or more KEEN waves interact, these edge structures would be affected first and therefore must be well resolved. 

As a strategy, we may choose a relatively large time step when exploring  a new regime of KEEN wave formation and  refine in velocity first and then in the space dimension until phase space error is under control. This can be assessed by the rate of global $L_2$ norm changes, for instance. Then, a further reduction in time step size which inevitably means interpolation errors will proliferate due to repeated applications of interpolations per physical time step. If the proper phase space structures are resolved, their evolution in time can be properly tracked with a small enough time step. But if the relevant scales are not being resolved, taking too small a time step will accelerate the degradation process.  This is why high order in time splitting is so crucial. It minimizes the number of iterations of interpolation required per unit gain in temporal accuracy. The largest time step possible for fixed error is the optimal strategy. But it must be found iteratively and only after proper phase space resolution which proceeds as fixed x, and v convergence inside the active region, followed by additional x refinement at fixed "converged" v resolution in the active region, followed by time step reduction till global measures of solution fidelity show little or no change. 

Note that the amplitude of the drive is multiplied by a factor $32$ between the canonical and the small drive cases. This can affect swings in the error in a similar way. Dividing the time step by two leads to a reduction of the error by a factor of  $2^6=64$ for a $6$-th order scheme. Such a reduction should therefore suffice. This argument implies that for a lower order in time scheme, much smaller time steps would be needed and as we mentioned, this would bring the error accumulation due to interpolation error to the fore. However, such arguments inherently rely on the degree of nonlinearity, or the number of scales needed to resolve in the sub-manifold to which the solution is confined. In terms of the smoothness of the solution, as afflicts all high order methods, simulations can not achieve full convergence when the solutions themselves are not smooth.  
  
\subsection{Convergence study via $\rho$-harmonics}
\begin{figure}[h!]
\begin{center}
\begin{tabular}{c}
\includegraphics[trim=1.cm 0.5cm 0cm 3cm, clip=true, width=0.5\linewidth,angle=0,height=5.5cm]{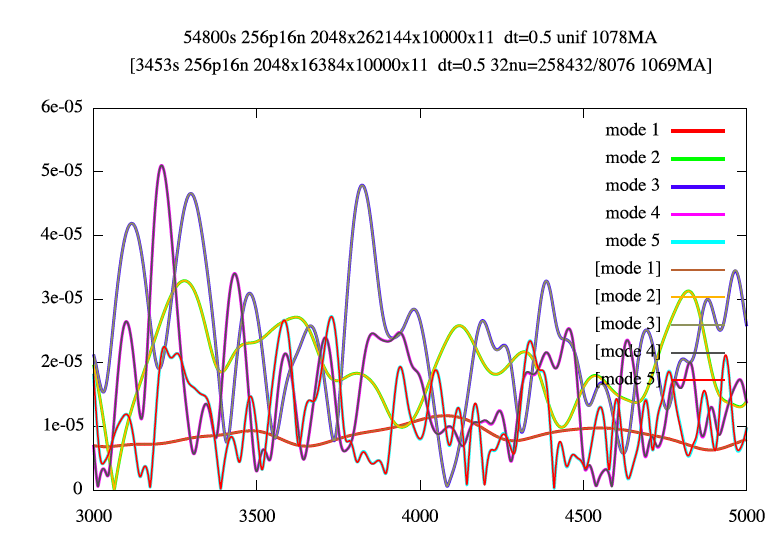} 
\end{tabular}
\end{center}
\caption{
Comparison of the time evolution of the first $5$ $\rho$ harmonics, in the {\bf small drive} amplitude case: run with {\bf uniform} velocity grid vs [canonical run with non-uniform velocity grid].
Parameters are $N_x=2048,\ N_v=262144,\ \Delta t=0.5$ uniform, $6th$ order time scheme (run12)
[$N_x=2048,\ N_v=16384,\ \Delta t=0.5$ non-uniform, $6$-th order time scheme (run13)].
}
\label{fig3_small_uniform}
\end{figure}

\begin{figure}[h!]
\begin{center}
\begin{tabular}{c}
\includegraphics[trim=1.cm 0.5cm 0cm 3cm, clip=true, width=0.45\linewidth,angle=0,height=5.cm]{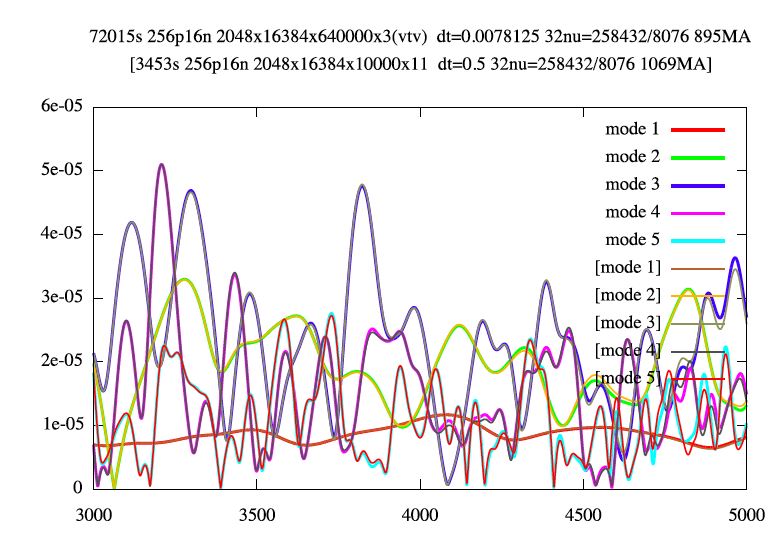} 
\end{tabular}
\end{center}
\caption{Comparison of the time evolution of the first $5$ $\rho$ harmonics, in the {\bf small drive} amplitude case: run with {\bf Strang} splitting scheme vs [canonical run with $6$th order time scheme].
Parameters are $N_x=2048,\ N_v=16384,\ \textrm{Strang},\  \Delta t=0.0078125$   non-uniform (run20)
[$N_x=2048,\ N_v=16384,\ \Delta t=0.5$ non-uniform, $6$-th order time scheme (run13)].}
\label{fig4_small_strang}
\end{figure}

\begin{figure}[h!]
\begin{center}
\begin{tabular}{c}
\includegraphics[trim=1.cm 0.5cm 0cm 3cm, clip=true, width=0.45\linewidth,angle=0,height=5.5cm]{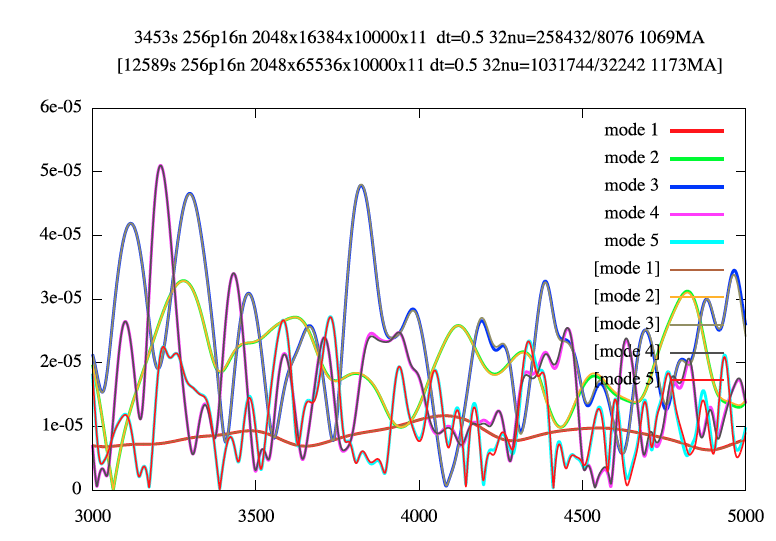}
\end{tabular}
\end{center}
\caption{Comparison of the time evolution of the first $5$ $\rho$ harmonics, in the {\bf small drive} amplitude case: canonical run vs. [{\bf velocity refined} run],
for testing the convergence in velocity. 
Parameters are
$N_x=2048,\ N_v=16384,\ \Delta t=0.5$ non-uniform, $6$-th order time scheme (run13) and 
[$N_x=2048,\ N_v=65536,\ \Delta t=0.5$ non-uniform, $6$-th order time scheme (run10)].
}
\label{fig5_small_refined}
\end{figure}

\begin{figure}[h!]
\begin{center}
\begin{tabular}{c}
\includegraphics[trim=1.cm 0.5cm 0cm 3cm, clip=true, width=0.45\linewidth,angle=0,height=5.cm]{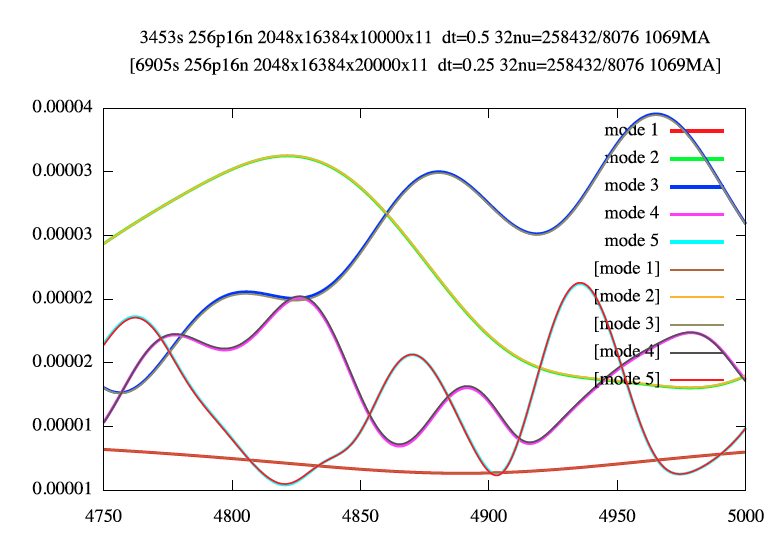}
\end{tabular}
\end{center}
\caption{Small drive:  $\Delta t=0.5$ vs [$\Delta t =0.25$]
Comparison of the time evolution of the first $5$ $\rho$ harmonics, in the {\bf small drive} amplitude case: canonical run vs [{\bf time refined} run], 
for testing the convergence in time. Parameters are
$N_x=2048,\ N_v=16384,\ \Delta t=0.5$ non-uniform,  $6$-th order time scheme (run13)
[$N_x=2048,\ N_v=16384,\ \Delta t=0.25$   non-uniform,  $6$-th order time scheme (run16a)].
}
\label{fig6_small_time_refined}
\end{figure}

\begin{figure}[h!]
\begin{center}
\begin{tabular}{c}
\includegraphics[trim=1.cm 0.5cm 0cm 3cm, clip=true, width=0.45\linewidth,angle=0,height=5.5cm]{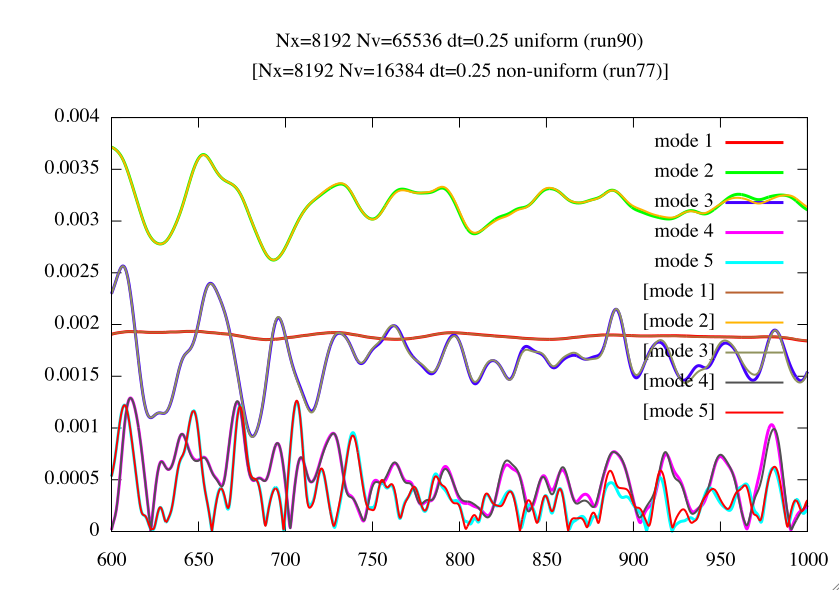}
\end{tabular}
\end{center}
\caption{Comparison of the time evolution of the first $5$ $\rho$ harmonics, in the {\bf canonical drive} case: run with {\bf uniform} velocity grid vs [canonical run with non-uniform velocity grid].
Parameters are $N_x=8192,\ N_v=65536,\ \Delta t=0.25$ uniform, $6th$ order time scheme (run90)
[$N_x=8192,\ N_v=16384,\ \Delta t=0.25$ non-uniform, $6$-th order time scheme (run77)].
}
\label{fig7_std_uniform}
\end{figure}

\begin{figure}[h!]
\begin{center}
\begin{tabular}{c}
\includegraphics[trim=1.cm 0.5cm 0cm 3cm, clip=true, width=0.45\linewidth,angle=0,height=5.5cm]{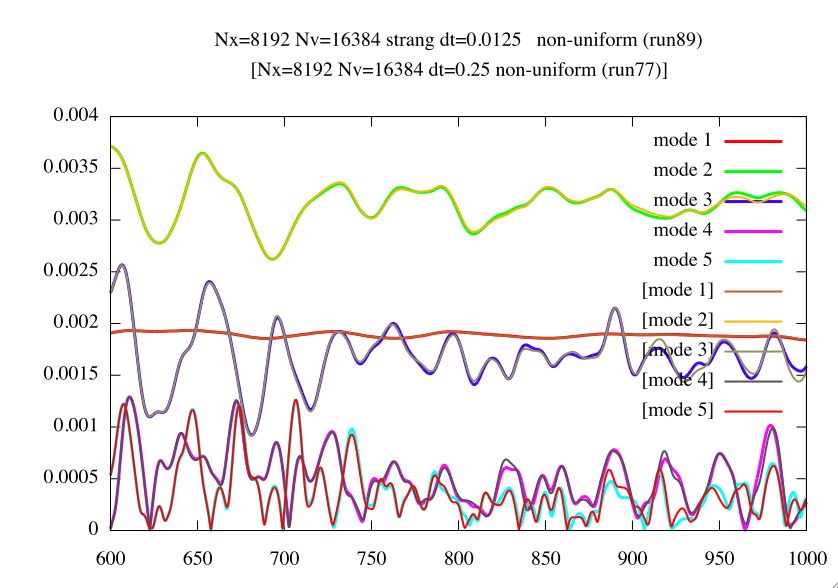}
\end{tabular}
\end{center}
\caption{Comparison of the time evolution of the first $5$ $\rho$ harmonics, in the {\bf canonical drive} amplitude case: run with {\bf Strang} splitting scheme vs [canonical run with $6$th order time scheme].
Parameters are $N_x=8192,\ N_v=16384,\ \textrm{Strang},\  \Delta t=0.0125$   non-uniform (run89)
[$N_x=8192,\ N_v=16384,\ \Delta t =0.25$ non-uniform, $6$-th order time scheme (run77)].}
\label{fig8_std_strang}
\end{figure}

We study the influence of numerical parameters, by comparing $\rho$-harmonics.
Results are given in Figures \ref{fig3_small_uniform} to \ref{fig6_small_time_refined}, for the small drive amplitude case,
and in Figures \ref{fig7_std_uniform} to \ref{fig10_std_time_refined}, for the canonical drive amplitude case. We would like to show that the density response, or spatial $\rho$-harmonics do not dependent on numerical parameters, so as to achieve, to the extent possible, converged results.  One general strategy here is to find a reasonably converged solution for a fixed time step, as shown in Figures \ref{fig1_small} and \ref{fig2_std} and to compare them with solutions deploying further refined grids in phase space or in time. All of this is done with the new method, that is non-uniform two grid cubic-splines interpolation in velocity, a fixed grid Lagrange interpolation of order $17$ in space and a $6$-th order symplectic time splitting integrator. On the other hand, we also compare the new method with more classical methods, such as uniform grids (with $6$-th order time splitting) and Strang splitting (with non-uniform grid). We choose the numerical parameters of the classical methods to be as small as possible leading to results of similar accuracy of the reasonably converged solution of the new method.

\subsubsection{The weakly driven case}
In Figure \ref{fig3_small_uniform}, we compare a uniform velocity grid run with $N_v=262144$ points which matches quite well with the reasonably converged solution ($N_v=16384$). 
But there is some discrepancy with respect to a refined run with the new method and $N_v=65536$ (see on Figure \ref{fig5_small_refined}). 
This can be explained by the fact that  $N_{\rm fine}=258432$, so that the mesh spacing of the fine grid is similar to the mesh spacing of the uniform grid with $N_v=262144$. In Figure \ref{fig4_small_strang}, we compare the Strang and $6$-th order schemes. The time step for Strang splitting is reduced to $\Delta t=0.0078125$ (we started with $\Delta t=0.25$ and then divided the time step successively by $2$ until we obtained an imperfect match). Note that taking so many time steps for the Strang splitting itself alters the solution. This could have been further checked by using a refined phase space solution; but at some point, one may begin to worry about the wall clock time of the simulation (here already at $20$ hours on $256$ processors). 
On the other hand, in Figure \ref{fig6_small_time_refined}, we see that results are not significantly affected by dividing the time step by $2$. This validates the fact that taking $\Delta t=0.5$ is good enough, when using the $6$-th order scheme. So, we gain a factor $262144/16384=16$ by using the non-uniform code, and a factor $(0.5/0.0078125)\cdot (2/11) \simeq 11.6$ by using the $6$-th order scheme
(we make the comparison here on the complexity and take $2$ steps for the Strang splitting, as the two $1/2$ step advections can be combined).
In order to avoid some undesirable effects that occur as a result of using too small a time step and a fine phase space grid, we should run a classical method with Strang splitting on a uniform grid, making the numerical resolution large enough to have similar accuracy. As a first estimate, a uniform velocity grid simulation with $N_v=262144$ points and  $\Delta t=0.0078125$ would cost
around $300$ hours (instead of $200\simeq 11.6\cdot 16\cdot 1\simeq 200$ hours [the reasonably converged simulation takes around $1$ hour only!] as we use $3$ steps for the Strang splitting, which eases the writing of the code, because visualization has to take place
on an integer multiple of time steps). In its present implementation,
that is $12$ days on $256$ processors; and we have limited the simulations here to require less than $24$ hours of continuous compute time (which can be done without the use of restart strategies).


\subsubsection{The canonical drive case}
\begin{figure}[h!]
\begin{center}
\begin{tabular}{c}
\includegraphics[trim=1.cm 0.5cm 0cm 3cm, clip=true, width=0.45\linewidth,angle=0,height=5.5cm]{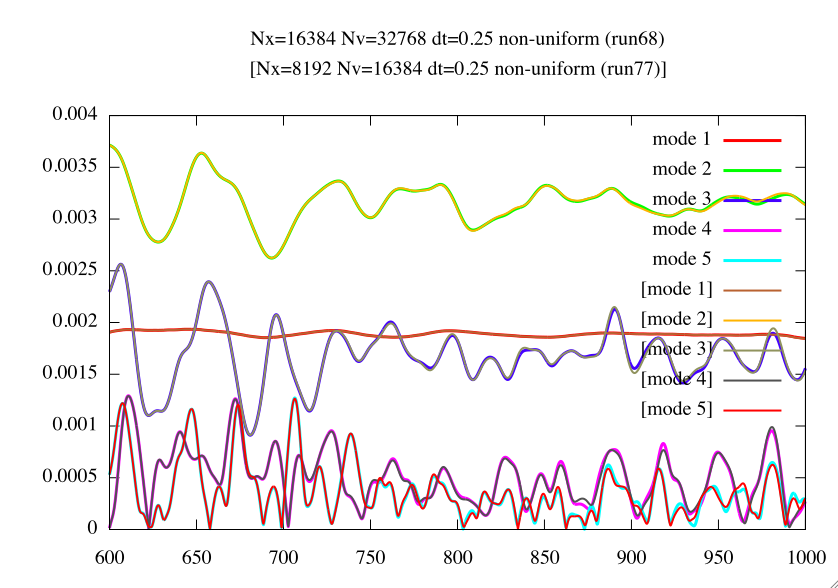}
\end{tabular}
\end{center}
\caption{Comparison of the time evolution of the first $5$ $\rho$ harmonics, in the {\bf canonical drive} amplitude case: {\bf phase-space refined} run vs [canonical run],
for testing the convergence in phase-space. 
Parameters are
$N_x=16384,\ N_v=32768,\ \Delta t=0.25$ non-uniform, $6$-th order time scheme (run68) and 
[$N_x=8192,\ N_v=16384,\  \Delta t=0.25$ non-uniform, $6$-th order time scheme (run77)].
}
\label{fig9_std_refined}
\end{figure}

\begin{figure}[h!]
\begin{center}
\begin{tabular}{c}
\includegraphics[trim=1.cm 0.5cm 0cm 3cm, clip=true, width=0.45\linewidth,angle=0,height=5.5cm]{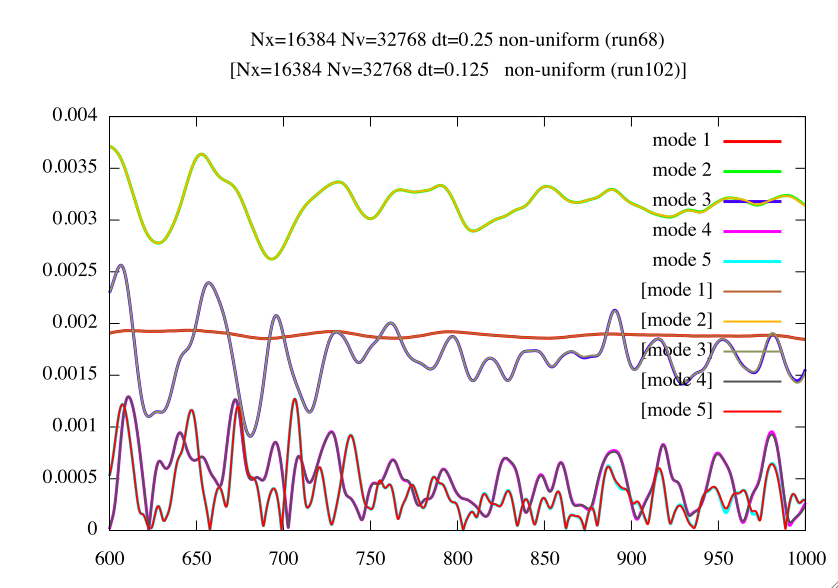}
\end{tabular}
\end{center}
\caption{Comparison of the time evolution of the first $5$ $\rho$ harmonics, in the {\bf canonical drive} amplitude case: phase-space refined run vs [phase-space and {\bf time refined} run], 
for testing the convergence in time. Parameters are
$N_x=16384,\ N_v=32768,\ \Delta t=0.25$ non-uniform,  $6$-th order time scheme (run68)
[$N_x=16384,\ N_v=32768,\ \Delta t=0.125$   non-uniform,  $6$-th order time scheme (run102)].}
\label{fig10_std_time_refined}
\end{figure}

Our investigations continued with the canonical drive case. Convergence of the reasonably refined run is checked in Figures  \ref{fig9_std_refined} and \ref{fig10_std_time_refined}.
More precisely, in Figure \ref{fig9_std_refined}, a comparison is made with a refined run using double the number of points in $x$ and $v$ and by taking the time step divided by $2$.
The time step is especially scrutinized in Figure \ref{fig10_std_time_refined}, using an already refined phase space simulation with $N_x=16384$ and $N_v=32768$. Note that the results match better in Figure \ref{fig10_std_time_refined} than in Figure  \ref{fig9_std_refined}, which indicates that lack of proper phase-space resolution seems to dominate time discretization error.

\noindent In Figure \ref{fig7_std_uniform}, comparison is made with a uniform run using $N_v=65536$ points. The match is not perfect; this may be due to the fact that
$N_{\rm fine}=89792$. On the other hand doubling the number of points in $v$ should lead to better results; note also that we are limited here to a grid size around
$2^{29}<8192\cdot 128172=2^{30}$ points, certainly because the distribution function is not stored in parallel. In Figure \ref{fig8_std_strang}, we compare with Strang splitting using 
$\Delta t =0.0125$; the result is even worse. We expect that we could take even a smaller time step in order to  achieve similar accuracy; but note that, we may have to refine the grids in phase space in order to prevent the accumulation of interpolation errors that may increase by taking more and more time steps; it is also interesting to see how the scheme behaves by showing results that are less converged.
 This point will be further studied with $L_2$-norm comparisons. The gain is here about $15$, which is less than $200$, but still valuable, and may be under-estimated. To conclude, we can say that we can reduce the computational costs by at least one order of magnitude for the canonical case, and by $2$ orders of magnitude for the small drive amplitude case.

\subsection{Time evolution of the $L_2$-norm}
\begin{figure}[h!]
\begin{center}
\begin{tabular}{c}
\includegraphics[trim=1.cm 0.5cm 0cm 0cm, clip=true, width=0.45\linewidth,angle=0,height=5.5cm]{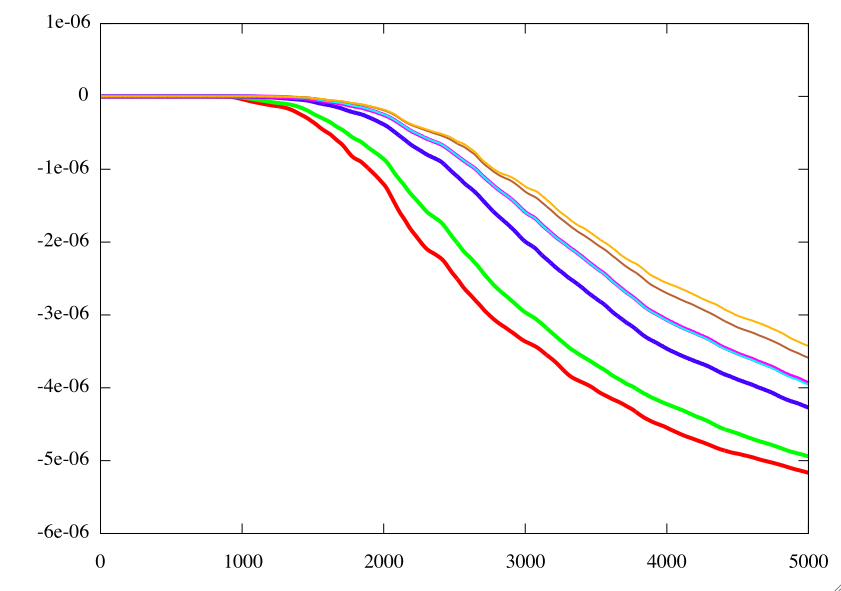}\\
\includegraphics[trim=0.cm 0.5cm 0cm 0cm, clip=true, width=0.45\linewidth,angle=0,height=3cm]{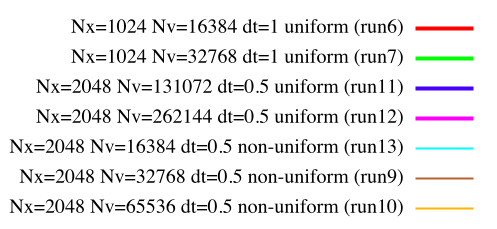}
\end{tabular}
\end{center}
\caption{Time evolution of relative $L_2$-norm for the {\bf small drive} amplitude case for different runs, with { uniform} and { non-uniform} velocity mesh.
The $6$-th order time splitting scheme is used for all the runs.
}
\label{fig11_l2_small}
\end{figure}

\begin{figure}[h!]
\begin{center}
\begin{tabular}{c}
\includegraphics[trim=1.cm 0.5cm 0cm 0cm, clip=true, width=0.45\linewidth,angle=0,height=5.5cm]{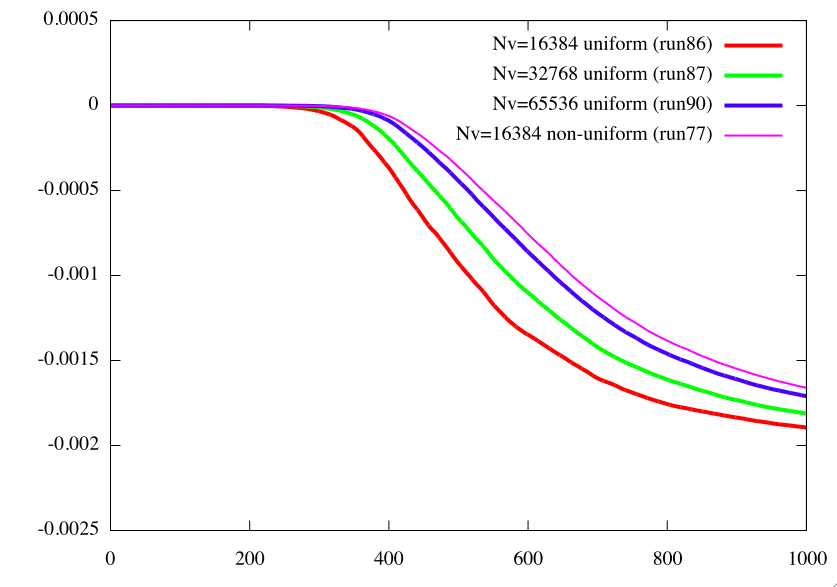} \\
\end{tabular}
\end{center}
\caption{Time evolution of relative $L_2$-norm for the canonical drive case: comparison of uniform grid runs with different velocity resolution
and a canonical non-uniform velocity grid run. Here $N_x=8192$ throughout.
}
\label{fig12_l2_std_uniform}
\end{figure}

\begin{figure}[h!]
\begin{center}
\begin{tabular}{c}
\includegraphics[trim=1.cm 0.5cm 0cm 0cm, clip=true, width=0.45\linewidth,angle=0,height=5.5cm]{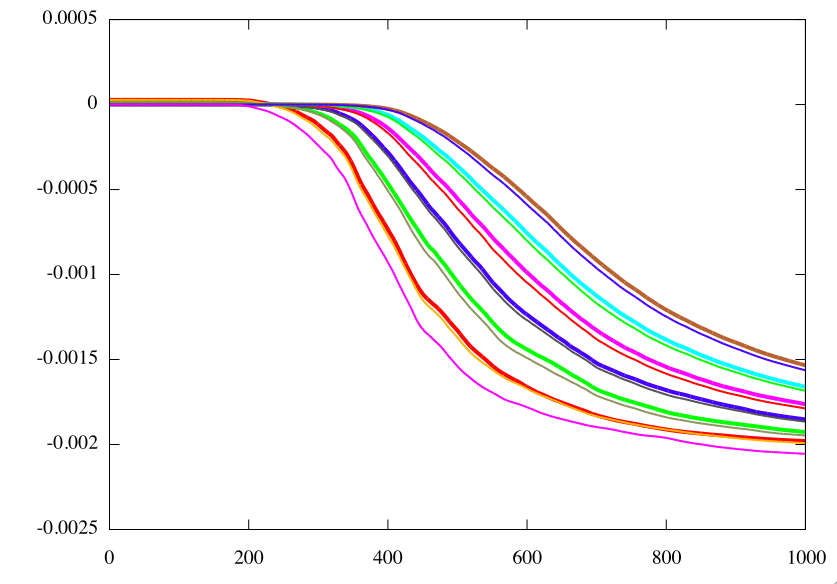} \\
\includegraphics[trim=1.cm 0.5cm 0cm 0cm, clip=true, width=0.45\linewidth,angle=0,height=4cm]{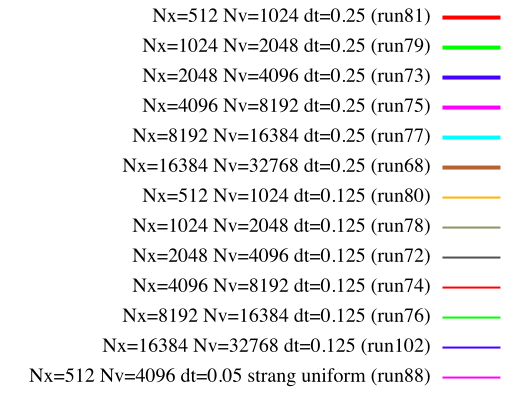} 
\end{tabular}
\end{center}
\caption{Time evolution of relative $L_2$-norm for the canonical drive amplitude case, with non-uniform velocity grids, $6$-th order time splitting and different phase-space resolutions
and time steps; only the last plot uses Strang splitting and uniform velocity grid.
}
\label{fig13_l2_std_refined}
\end{figure}

\begin{figure}[h!]
\begin{center}
\begin{tabular}{c}
\includegraphics[trim=1.cm 0.5cm 0cm 0cm, clip=true, width=0.45\linewidth,angle=0,height=5.5cm]{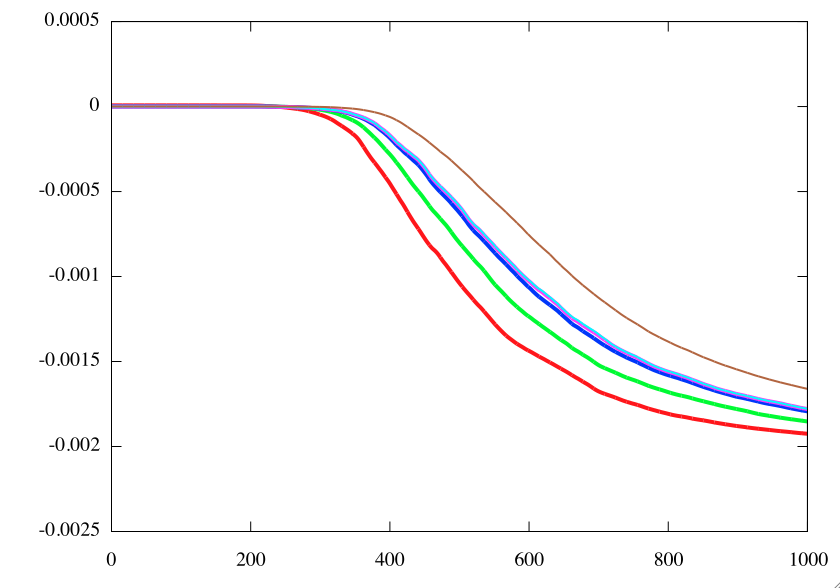} \\
\includegraphics[trim=0.cm 0.cm 0cm 0cm, clip=true, width=0.2\linewidth,angle=0,height=2.cm]{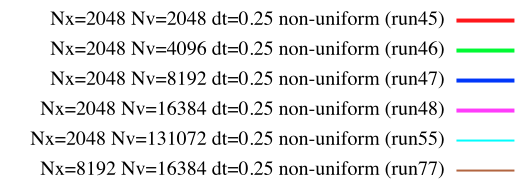}
\end{tabular}
\end{center}
\caption{Time evolution of relative $L_2$-norm for the canonical drive amplitude case: comparison of non-uniform runs with different velocity resolution
and {\bf space resolution} {\bf $N_x=2048$}, and with canonical non-uniform run, using $N_x=8192$. Here $\Delta t =0.25$.
}
\label{fig14_l2_std_space}
\end{figure}

We study the influence of  the $L_2$-norm conservation on the numerical parameters.
Results are given on Figure \ref{fig11_l2_small} for the small drive amplitude case. A more detailed study is performed on the 
canonical drive amplitude case. Results are shown in Figures \ref{fig12_l2_std_uniform} to \ref{fig17_l2_std_6th_zoom}.
The study of the conservation of $L_2$-norm gives an indication of the global error that is made. It is clearly not enough as some numerical schemes can conserve the $L_2$-norm exactly while otherwise proliferating of errors.
However the lower bound of the $L_2$ error has the advantage of not needing a converged solution. 
So, we emphasize that this study is quite instructive, as it has permitted us to detect that taking $N_x=2048$ points in $x$ was not enough at a certain point, 
and we thank the anonymous referee for suggesting us to add such a diagnostic. Error computation on converged solutions should be the next complementary step for checking the convergence. Such a step may need more elaborate post-processing tools as solution may differ with the number of points in $x$, $v$ and time step; choices have to be made, and this may not be obvious leading to certain proliferation of work especially for computational studies producing  large data sets. 

\subsubsection{The weakly driven case}

\noindent In Figure \ref{fig11_l2_small}, we compare uniform and non-uniform velocity grid runs and change the number of points in velocity, using always the $6$-th order time splitting scheme.
Time step is fixed at $0.5$ (only for some uniform runs with poorest resolution we choose a time step twice as big). We clearly see that in order to have as good an $L_2$-norm
we have to take $N_v$ $16$ times bigger in the uniform case. Taking $N_v$ even bigger, that is $N_v=65536,$ improves the $L_2$-norm in the non-uniform run. Note that this would correspond
to a uniform mesh with $N_{\rm fine}=1031744$ points. We can however note that the gain in $L_2$-norm is then not so big, when compared to non-uniform runs with $N_v=32768$.
Certainly, we should deploy also more points in $x$ to improve the results, as we will see on the more detailed study for the canonical drive case.

\subsubsection{The canonical drive case}
\noindent In Figure \ref{fig12_l2_std_uniform}, we compare the $L_2$-norm of the reasonably converged non-uniform run with uniform runs. We again observe that taking
$N_v=65536$ for the uniform run is less accurate than taking $N_v=16384$ points for the non-uniform run. We also note that the $L_2$-norm is less precisely conserved as in the small drive case, which
again confirms that it is more challenging to have converged results for large amplitude KEEN waves which shed small vortices and generally pursue higher order dynamics that the slower evolving and smaller amplitude response of the plasma when driven much more weakly.  

\noindent In Figure \ref{fig13_l2_std_refined}, we compare different non-uniform runs, with $\Delta t=0.25$ and $\Delta t=0.125$, using the $6$-th order time splitting scheme.
This run is compared with a classical run (Strang splitting and a uniform mesh), with classical parameters ($\Delta t=0.05,\ N_x=512,\ N_v=4096$). We clearly see the improvements
of the non-uniform runs. 
The level of $L_2$-norm accuracy  achieved  at time $T=400$ for the classical run, is maintained at time $T\simeq700$ 
 for the reasonable converged run and at time $T\simeq800$ for the refined run.
 We also note that dividing the time step by $2$ leads to a degradation of the $L_2$-norm  but only a small one; a further study on the choice of proper time step will be discussed later in the paper.

\noindent In Figure \ref{fig14_l2_std_space}, we plot the $L_2$-norm for different runs with fixed spatial resolution $N_x=2048$. We clearly see that the $L_2$-norm stucks at a certain level, and that changing the number of points in $v$ has no more real effect, when $N_v \ge 16384$ (or even $N_v=8192$): taking $N_v=131072$ does not lead to a real improvement. 
On the contrary, taking more points in $x$ has a clear effect, as we can see on the run using $N_x=8192$ and $N_v=16384$.

\begin{figure}[h!]
\begin{center}
\begin{tabular}{c}
\includegraphics[trim=1.cm 0.5cm 0cm 0cm, clip=true, width=0.45\linewidth,angle=0,height=5.5cm]{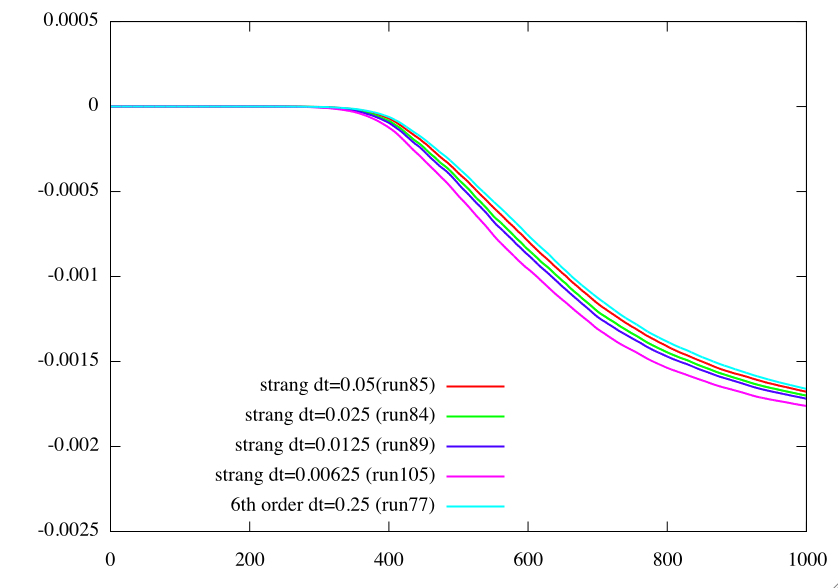}
\end{tabular}
\end{center}
\caption{Time evolution of relative $L_2$-norm for the canonical drive amplitude case: comparison of non-uniform runs with {\bf Strang} splitting and different time steps; and also with canonical non-uniform run with $6$th order scheme. Here $N_x=8192,\ N_v=16384$.
}
\label{fig15_l2_std_strang}
\end{figure}

\begin{figure}[h!]
\begin{center}
\begin{tabular}{c}
\includegraphics[trim=1.cm 0.5cm 0cm 0cm, clip=true, width=0.45\linewidth,angle=0,height=5.5cm]{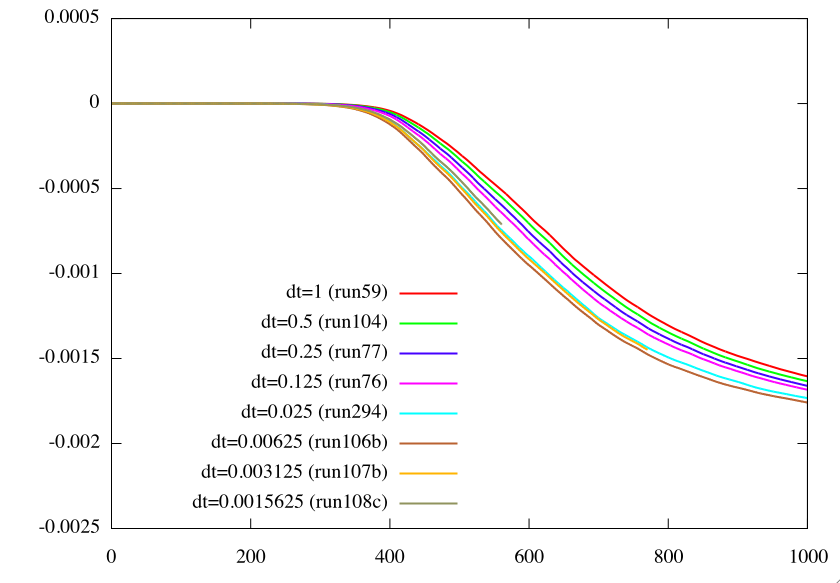}
\end{tabular}
\end{center}
\caption{Time evolution of relative $L_2$-norm for the canonical drive amplitude case: comparison of non-uniform runs with {\bf $6$th order} scheme and different time steps.
Here $N_x=8192,\ N_v=16384$.
}
\label{fig16_l2_std_6th}
\end{figure}

\begin{figure}[h!]
\begin{center}
\begin{tabular}{c}
\includegraphics[trim=1.cm 0.5cm 0cm 0cm, clip=true, width=0.45\linewidth,angle=0,height=5.5cm]{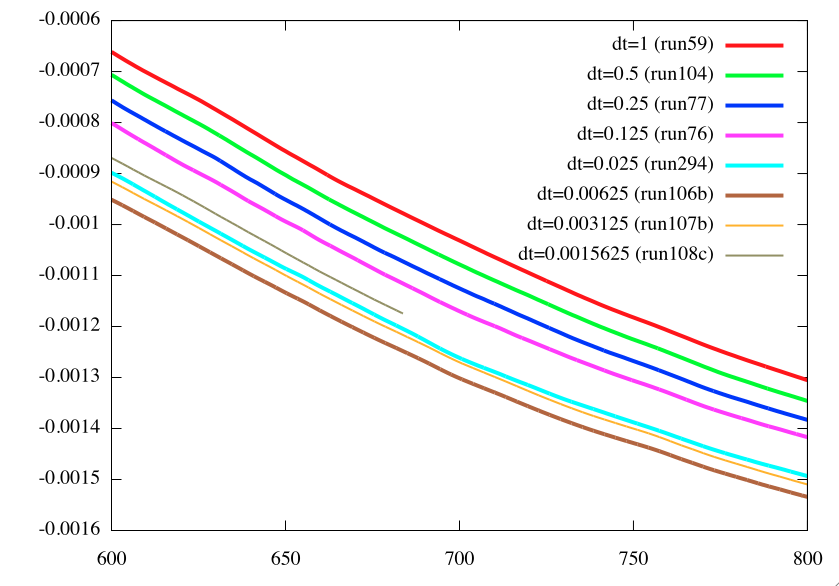}
\end{tabular}
\end{center}
\caption{Zoom of previous picture for $t\in [600,800]$.
}
\label{fig17_l2_std_6th_zoom}
\end{figure}

\noindent The last plots (Figures \ref{fig15_l2_std_strang} to \ref{fig17_l2_std_6th_zoom}) concern the study of the influence of the time step. 

\noindent We first consider on Figure  \ref{fig15_l2_std_strang}, different runs using the Strang splitting scheme. We see that decreasing the time step leads to a decrease of the $L^2$-norm conservation. Note that the situation is even worse between $\Delta t=0.00625$ and $\Delta t=0.0125$, 
than between $\Delta t=0.025$ and $\Delta t=0.0125$. This may be explained by the fact that we use more and more interpolations and that the underlying phase-space resolution is not enough
(here $N_x=8192,\ N_v=16384$). Note that using the $6$-th order scheme with $\Delta t=0.25$ leads to better $L_2$-norm conservation as less time iterations are needed.
In other words, Strang splitting scheme needs here smaller time steps to converge in time; but on the other hand, taking smaller time step leads to more interpolation and thus
more errors, when the grid is not refined enough. 

\noindent In Figures \ref{fig16_l2_std_6th} and \ref{fig17_l2_std_6th_zoom}, we consider the $6$-th order scheme and look for the influence of the time step. It is interesting to see that 
the $L_2$-norm decreases with the time step and that after a while (taking time step smaller than on the previous plot with Strang splitting), 
the $L_2$-norm begins to increase again. One explanation may come from the fact that the time step is so small that the foot of the characteristic is very near to the grid point, so that the interpolation error
is reduced, or at least the $L_2$-norm is better conserved (see \cite{ChaDeMe2012} for details on such phenomenon), depending on the type of interpolation that is used.

Note that cubic splines lead to exact $L^2$-norm conservation at the limit of a $0$ displacement (exponential integrator), but induce numerical dispersion, just as centered reconstructions in finite volume schemes.
On the other hand, odd order Lagrange interpolation (as is used here in the $x$ direction) acts just as an upwind scheme, which gives better stability, when using small displacements, without tending towards  an exact $L_2$-norm conservation in the $0$ displacement limit. Further studies, such as detailed quantitative error estimates when compared to converged solutions should confirm that there is an optimal time step which is near $\Delta t = 0.25$ for the problems treated herein. It is not possible to reach such a conclusion by merely looking at the global $L_2$-norm scaling with time.


\subsection{Weakly driven cases with different drive times.}
\begin{figure}[h!]
\begin{tabular}{c}
\includegraphics[trim=0.5cm 4cm 2cm 2.4cm, clip=true, width=0.45\linewidth,angle=0,height=5.5cm]{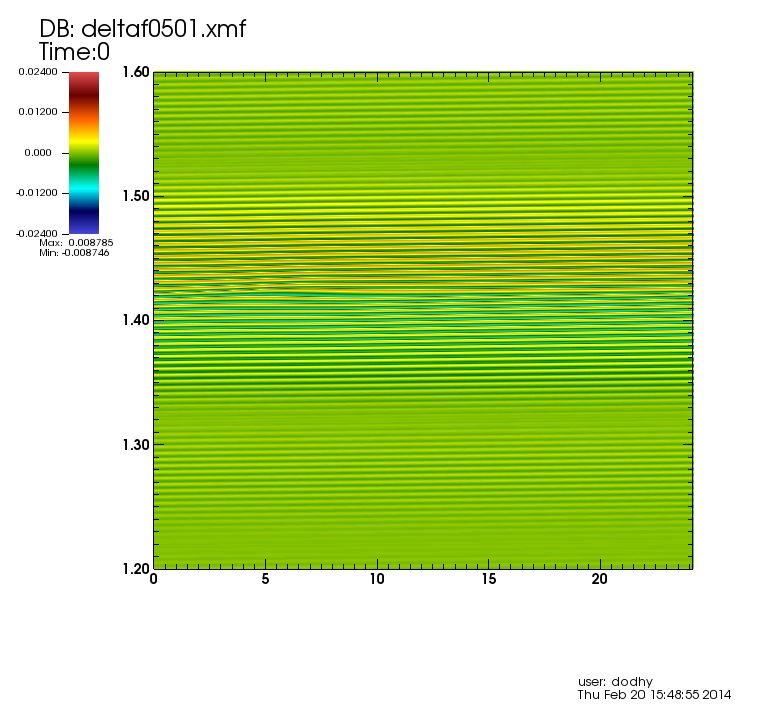}\\
\includegraphics[trim=0.5cm 4cm 2cm 2.4cm, clip=true, width=0.45\linewidth,angle=0,height=5.5cm]{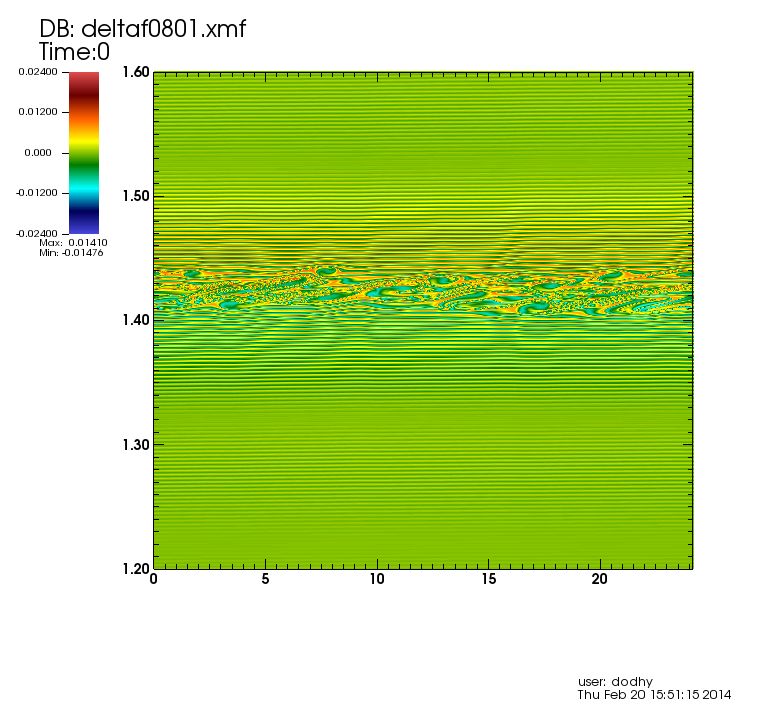}\\
\includegraphics[trim=0.5cm 4cm 2cm 2.4cm, clip=true, width=0.45\linewidth,angle=0,height=5.5cm]{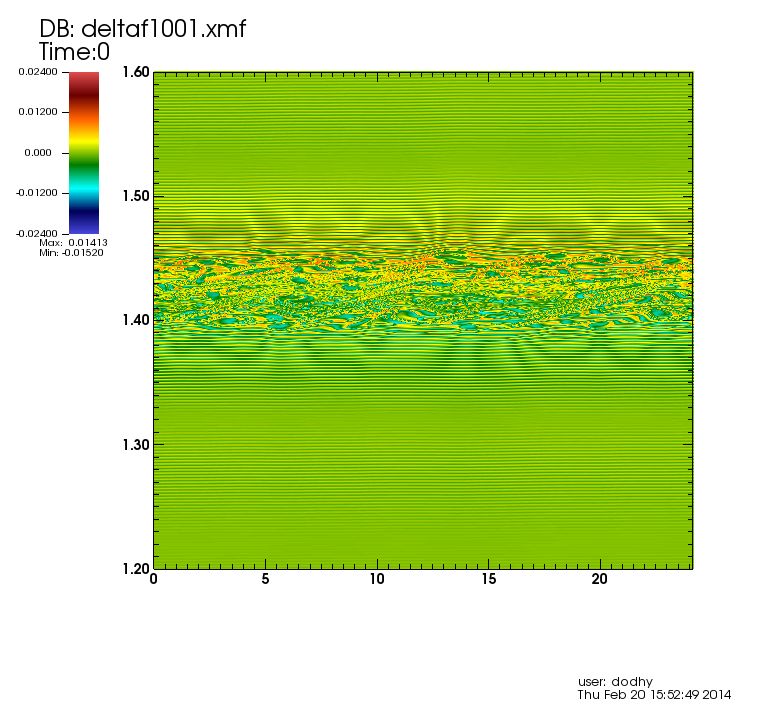}\\
\includegraphics[trim=1.cm 0.5cm 0cm 2.5cm, clip=true, width=0.45\linewidth,angle=0,height=5cm]{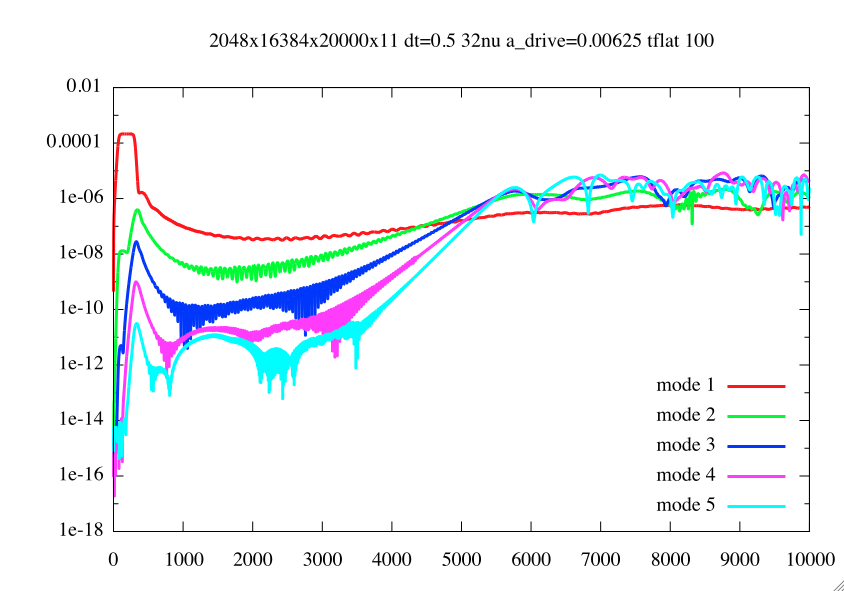} 
\end{tabular}
\caption{Case $a_{\rm Dr}=0.00625$ and $T_{\rm Dr}=100$.
$\delta f$ distribution function $f-f_0$ at time $T=5000,\ 8000,\ 10000$ as a function of $(x,v) \in [0,4\pi]\times [1.2,1.6]$ (from top to bottom) and time evolution of the first $5$ $\rho$ harmonics (bottom): 
canonical run with $6$th order time scheme and non-uniform velocity mesh.
Parameters are $N_x=2048,\ N_v=16384,\ \Delta t=0.5$ non-uniform, $6$-th order time scheme (run13).}
\label{Tdrive100}
\end{figure}

\begin{figure}[h!]
\begin{tabular}{c}
\includegraphics[trim=0.5cm 4cm 2cm 2.4cm, clip=true, width=0.45\linewidth,angle=0,height=5.5cm]{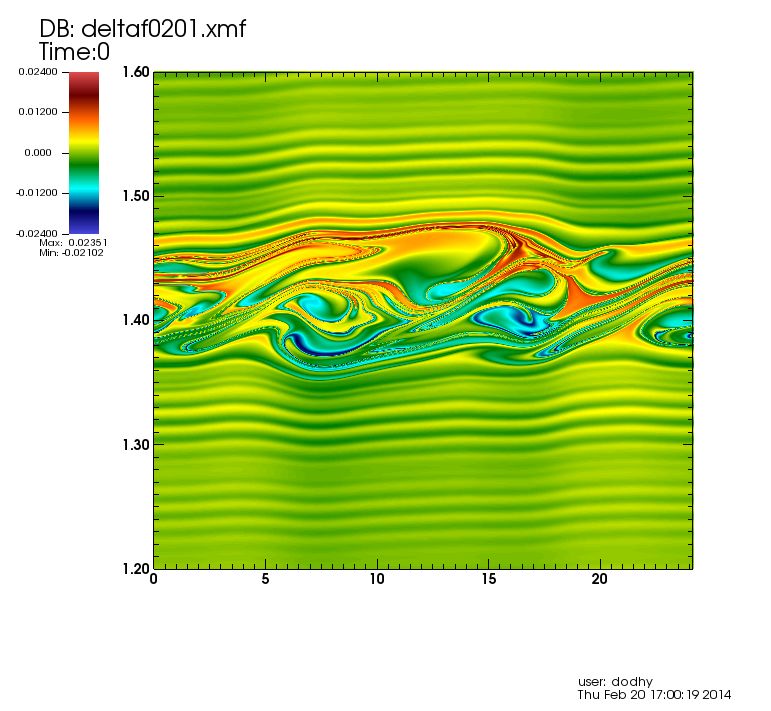}\\
\includegraphics[trim=0.5cm 4cm 2cm 2.4cm, clip=true, width=0.45\linewidth,angle=0,height=5.5cm]{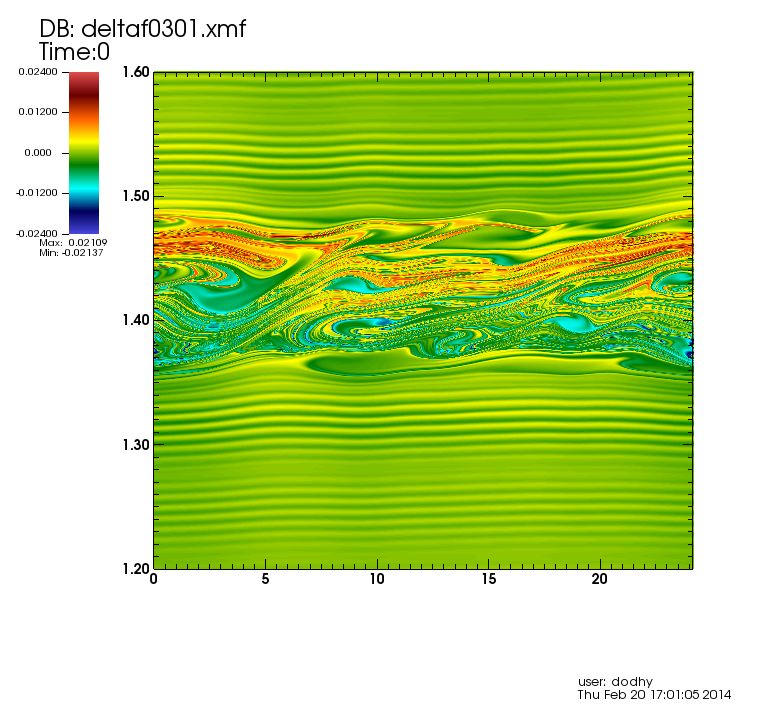}\\
\includegraphics[trim=0.5cm 4cm 2cm 2.4cm, clip=true, width=0.45\linewidth,angle=0,height=5.5cm]{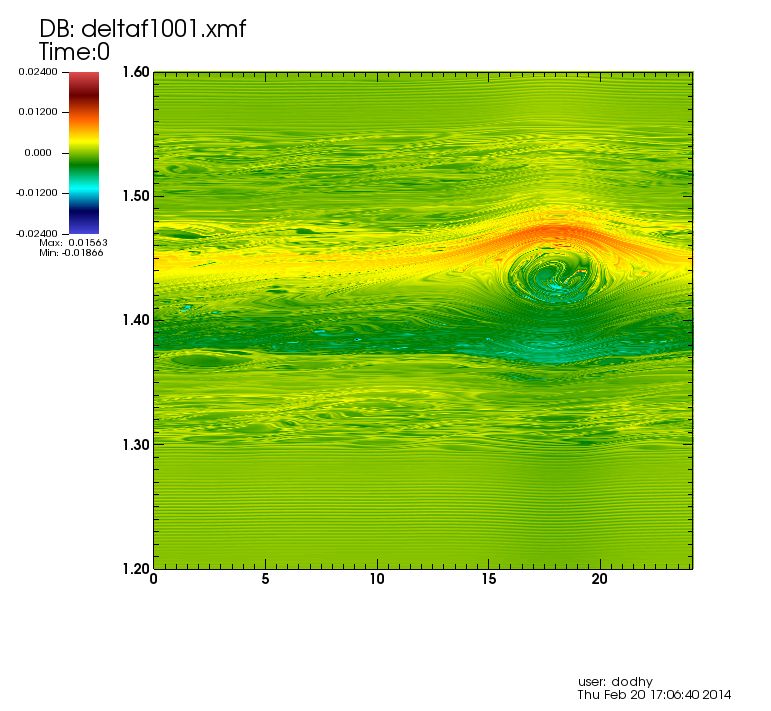}\\
\includegraphics[trim=1.cm 0.5cm 0cm 2.5cm, clip=true, width=0.45\linewidth,angle=0,height=5.cm]{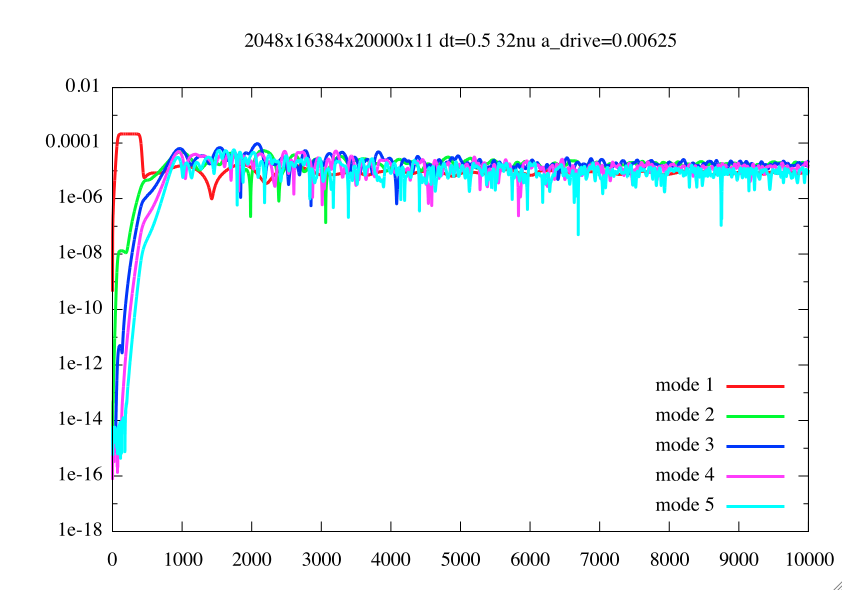} 
\end{tabular}
\caption{Case $a_{\rm Dr}=0.00625$ and $T_{\rm Dr}=200$.
$\delta f$ distribution function $f-f_0$ at time $T=2000,\ 3000,\ 10000$ as a function of $(x,v) \in [0,4\pi]\times [1.2,1.6]$ (from top to bottom; from left to right) and time evolution of the first $5$ $\rho$ harmonics (bottom): 
canonical run with $6$th order time scheme and non-uniform velocity mesh.
Parameters are $N_x=2048,\ N_v=16384,\ \Delta t=0.5$ non-uniform, $6$-th order time scheme (run13).}
\label{Tdrive200}
\end{figure}

\begin{figure}[h!]
\begin{tabular}{c}
\includegraphics[trim=0.5cm 4cm 2cm 2.4cm, clip=true, width=0.25\linewidth,angle=0,height=5.5cm]{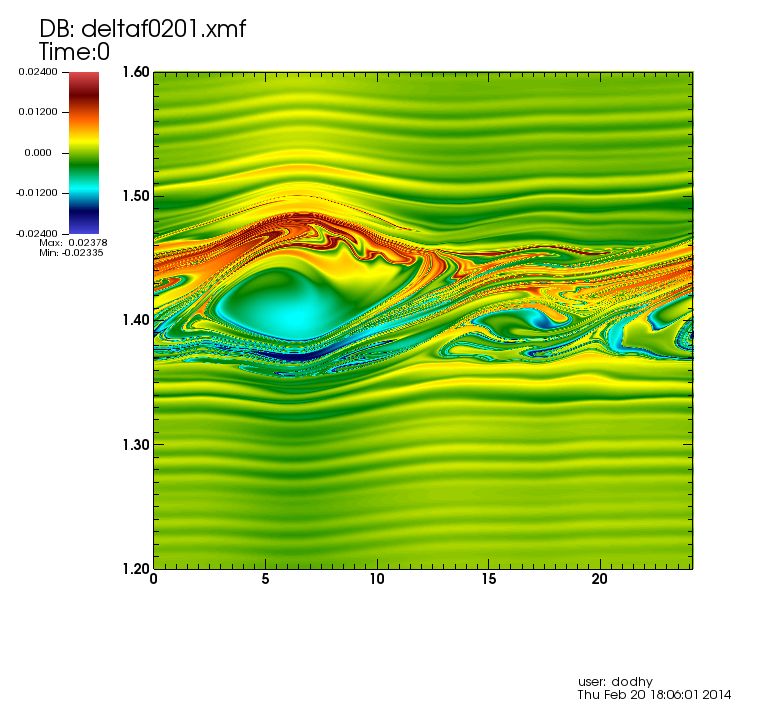}\\
\includegraphics[trim=0.5cm 4cm 2cm 2.4cm, clip=true, width=0.25\linewidth,angle=0,height=5.5cm]{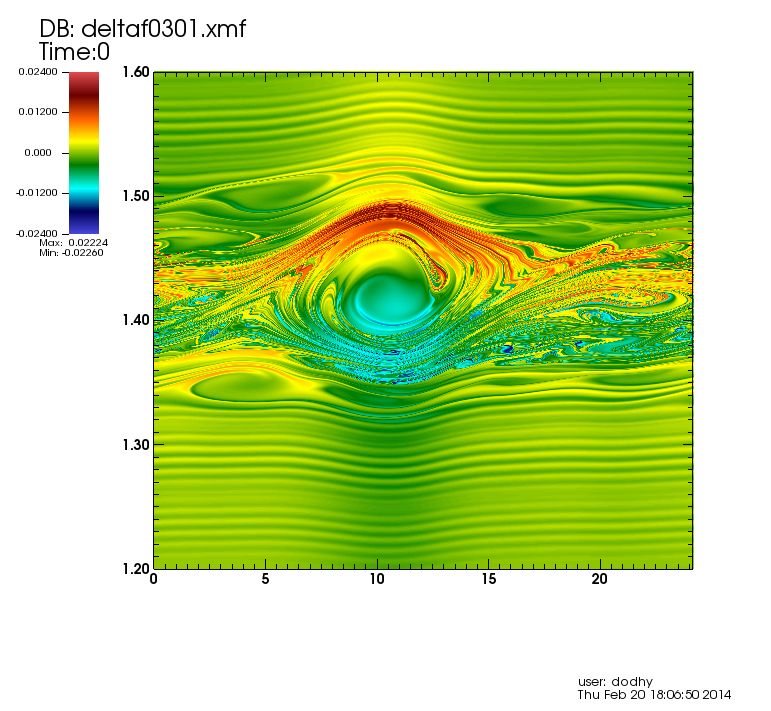}\\
\includegraphics[trim=0.5cm 4cm 2cm 2.4cm, clip=true, width=0.25\linewidth,angle=0,height=5.5cm]{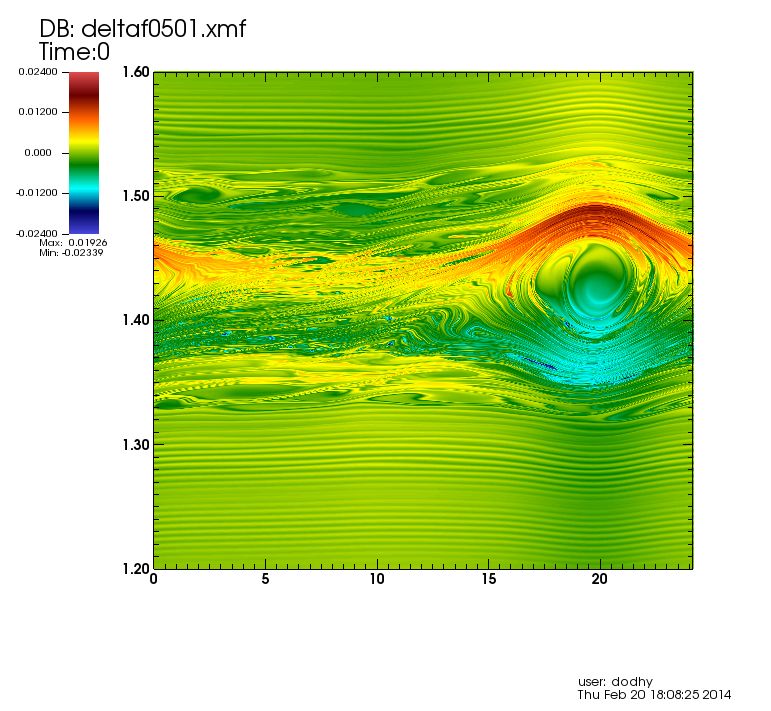}\\
\includegraphics[trim=1.cm 0.5cm 0cm 2.5cm, clip=true,  width=0.45\linewidth,angle=0,height=5cm]{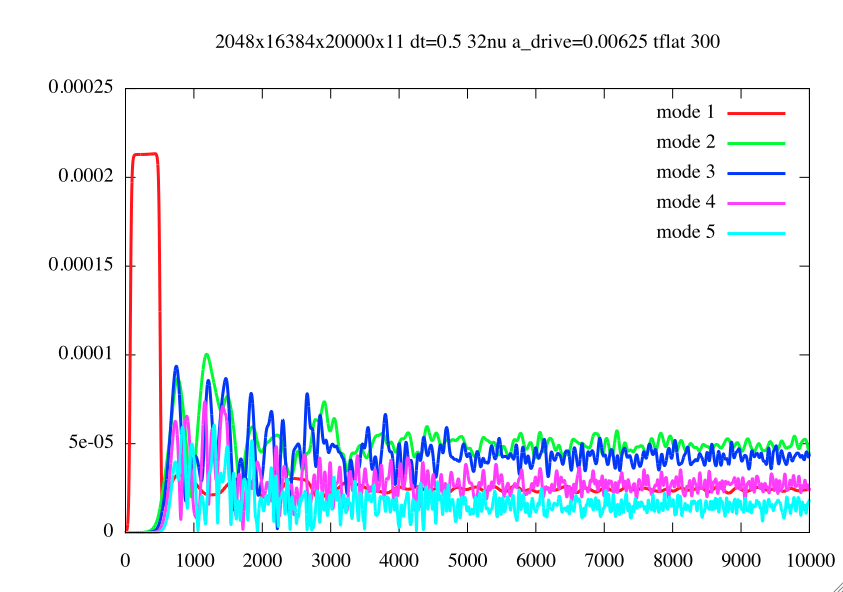} 
\end{tabular}
\caption{Case $a_{\rm Dr}=0.00625$ and $T_{\rm Dr}=300$.
$\delta f$ distribution function $f-f_0$ at time $T=2000,\ 3000,\ 5000$ as a function of $(x,v) \in [0,4\pi]\times [1.2,1.6]$ (from top to bottom; from left to right) and time evolution of the first $5$ $\rho$ harmonics (bottom): 
canonical run with $6$th order time scheme and non-uniform velocity mesh.
Parameters are $N_x=2048,\ N_v=16384,\ \Delta t=0.5$ non-uniform, $6$-th order time scheme (run33).}
\label{Tdrive300}
\end{figure}

\begin{figure}[h!]
\begin{tabular}{c}
\includegraphics[trim=0.cm 0cm 0cm 0cm, clip=true, width=0.25\linewidth,angle=0,height=5cm]{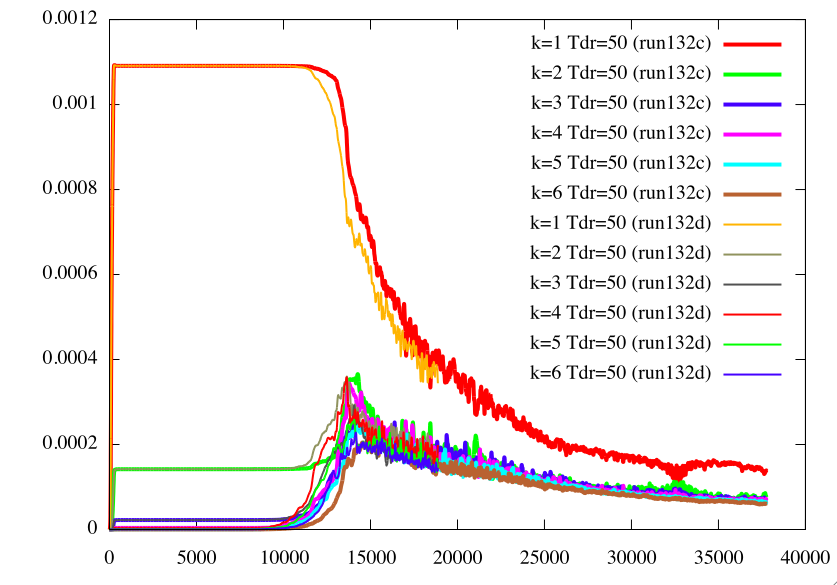}\\
\includegraphics[trim=0.cm 0cm 0cm 0cm, clip=true, width=0.25\linewidth,angle=0,height=5cm]{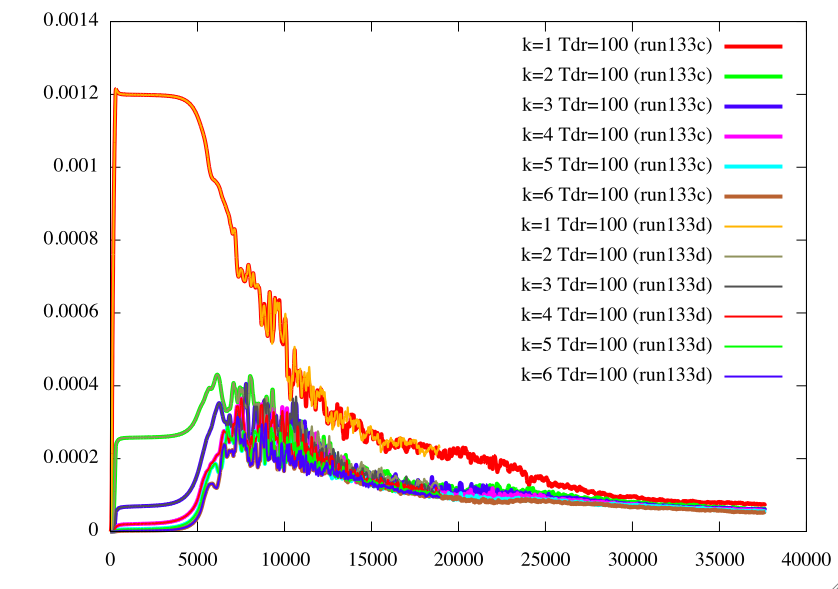}\\
\includegraphics[trim=0.cm 0cm 0cm 0cm, clip=true, width=0.25\linewidth,angle=0,height=5cm]{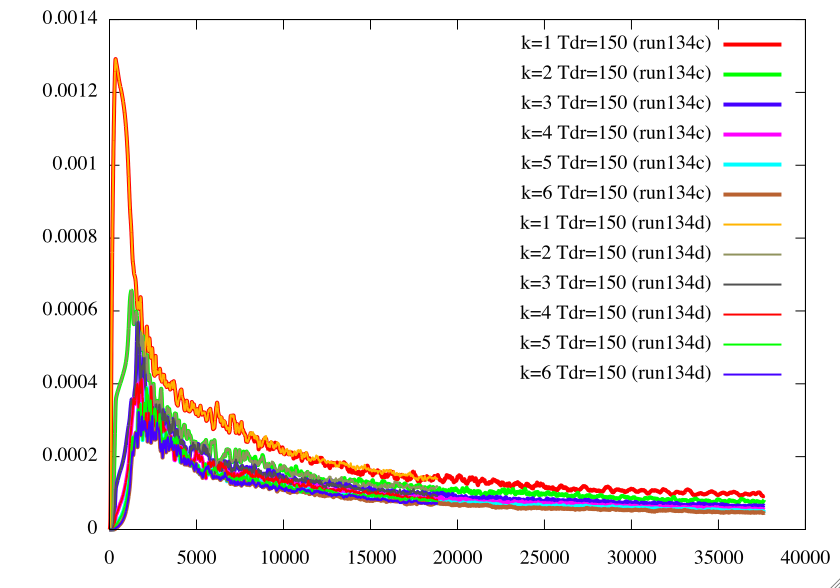}\\
\includegraphics[trim=0.cm 0cm 0cm 0cm, clip=true, width=0.25\linewidth,angle=0,height=5cm]{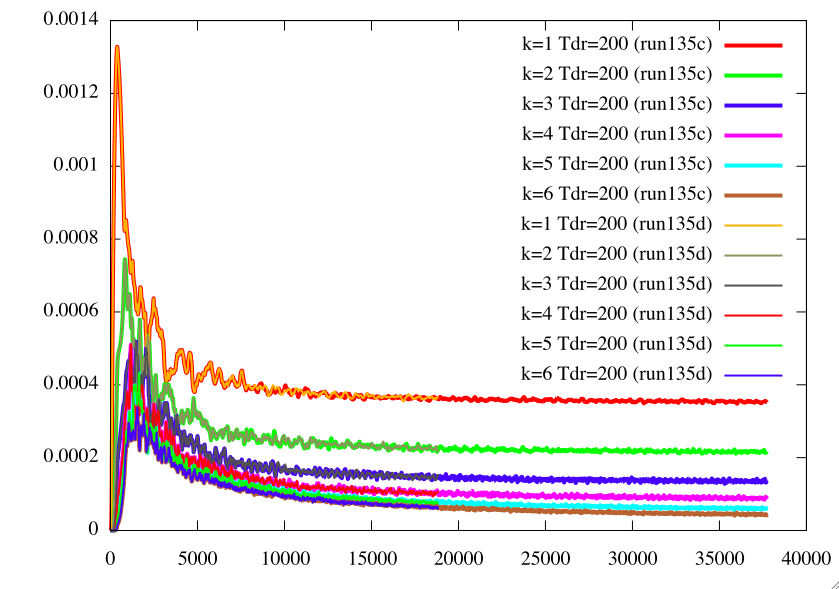}
\end{tabular}
\caption{Time evolution of $\sqrt{\int_\mathbb{R}|\hat{f}^{(k)}(t,v)|^2dv}$ for $k=1,\ 2,\dots, 6$ and for drive time $50,\ 100,\ 150,\ 200$ (run132c-135c and run132d-135d)
Parameters are $N_x=4096,\ N_v=32768$, non-uniform, $6$-th order time scheme and $\Delta t=0.5$ for run132c-135c, $\Delta t=0.25$ for run132d-135d.}
\label{fig21}
\end{figure}

\begin{figure}[h!]
\begin{tabular}{c}
\includegraphics[trim=0.cm 0cm 0cm 0cm, clip=true, width=0.25\linewidth,angle=0,height=5cm]{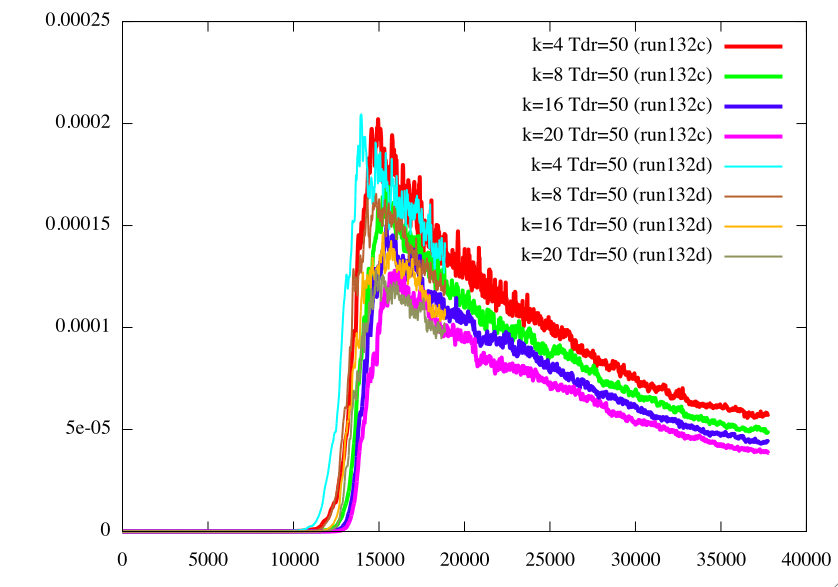}\\
\includegraphics[trim=0.cm 0cm 0cm 0cm, clip=true, width=0.25\linewidth,angle=0,height=5cm]{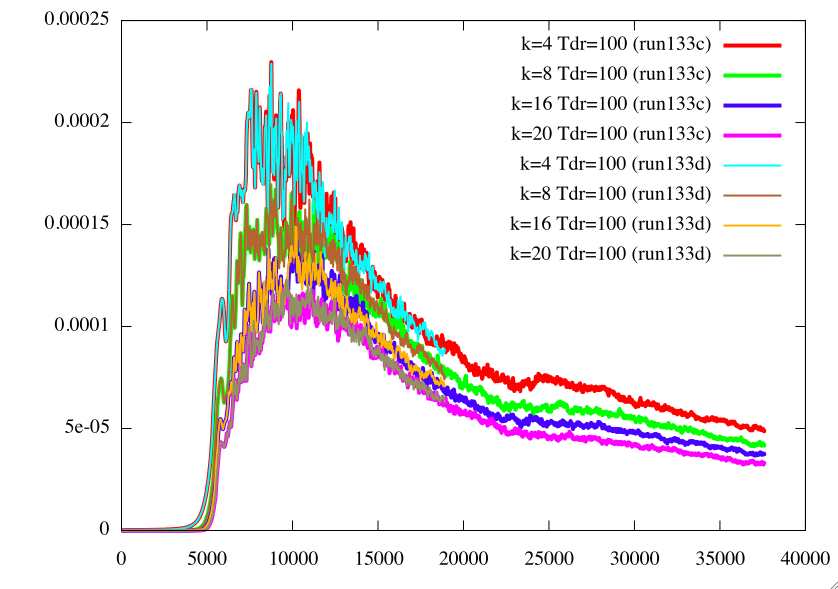}\\
\includegraphics[trim=0.cm 0cm 0cm 0cm, clip=true, width=0.25\linewidth,angle=0,height=5cm]{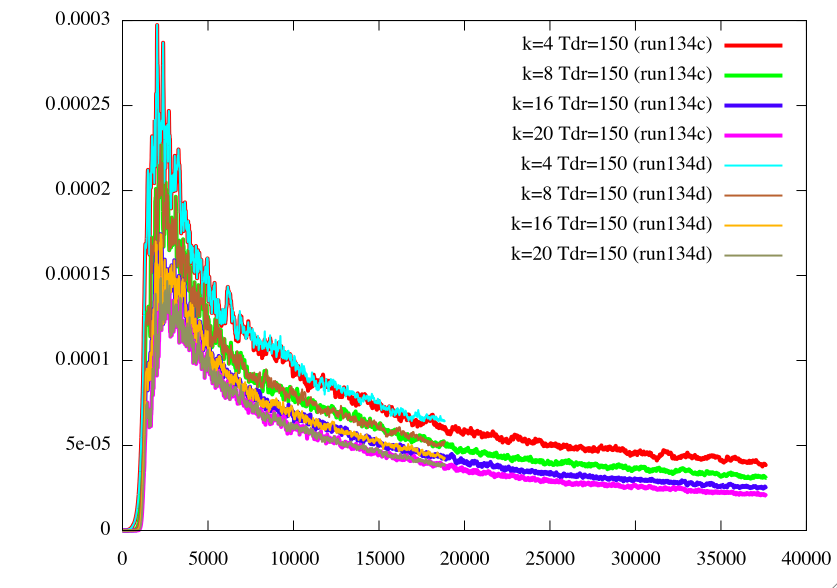}\\
\includegraphics[trim=0.cm 0cm 0cm 0cm, clip=true, width=0.25\linewidth,angle=0,height=5cm]{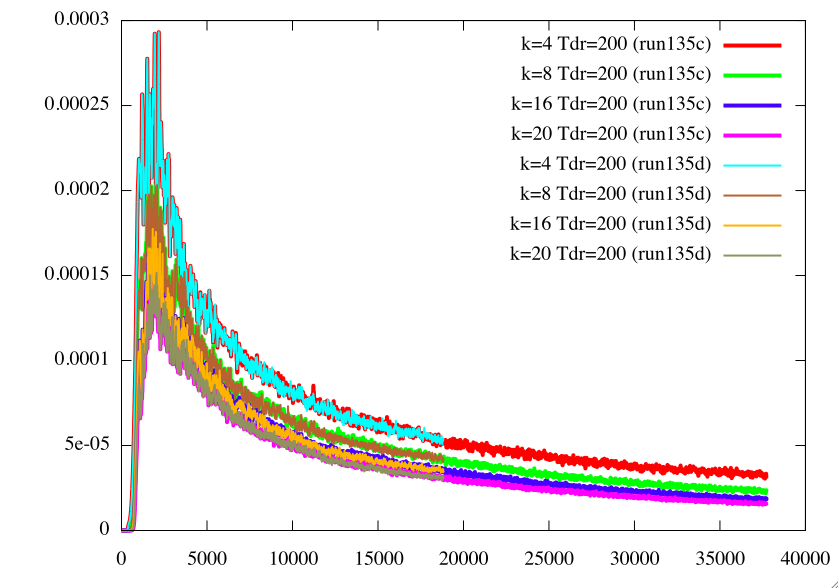}
\end{tabular}
\caption{Time evolution of $\sqrt{\int_\mathbb{R}|\hat{f}^{(k)}(t,v)|^2dv}$ for $k=4,\ 8,\ 16,\ 20$ and for drive time $50,\ 100,\ 150,\ 200$ (run132c-135c and run132d-135d)
Parameters are $N_x=4096,\ N_v=32768$, non-uniform, $6$-th order time scheme and $\Delta t=0.5$ for run132c-135c, $\Delta t=0.25$ for run132d-135d.}
\label{fig22}
\end{figure}

\begin{figure}[h!]
\begin{tabular}{c}
\includegraphics[trim=0.cm 0cm 0cm 0cm, clip=true, width=0.25\linewidth,angle=0,height=5cm]{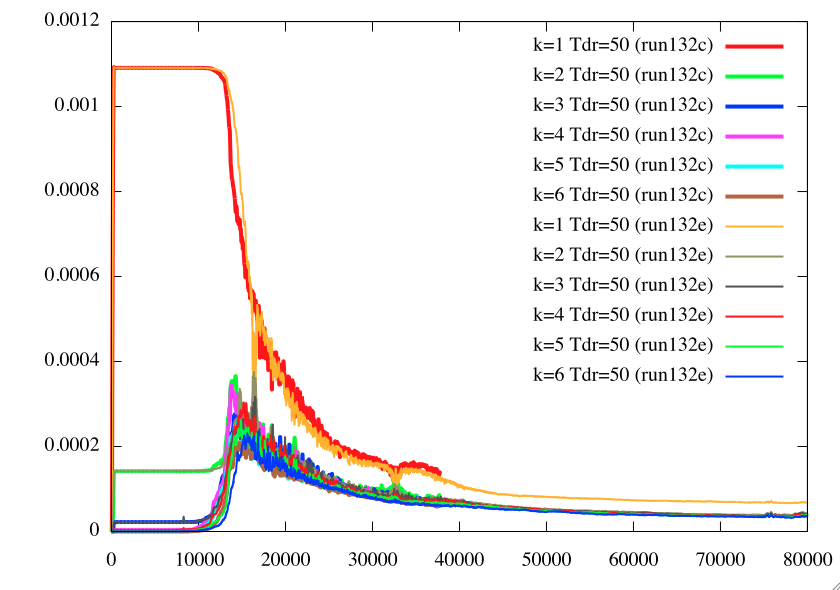}\\
\includegraphics[trim=0.cm 0cm 0cm 0cm, clip=true, width=0.25\linewidth,angle=0,height=5cm]{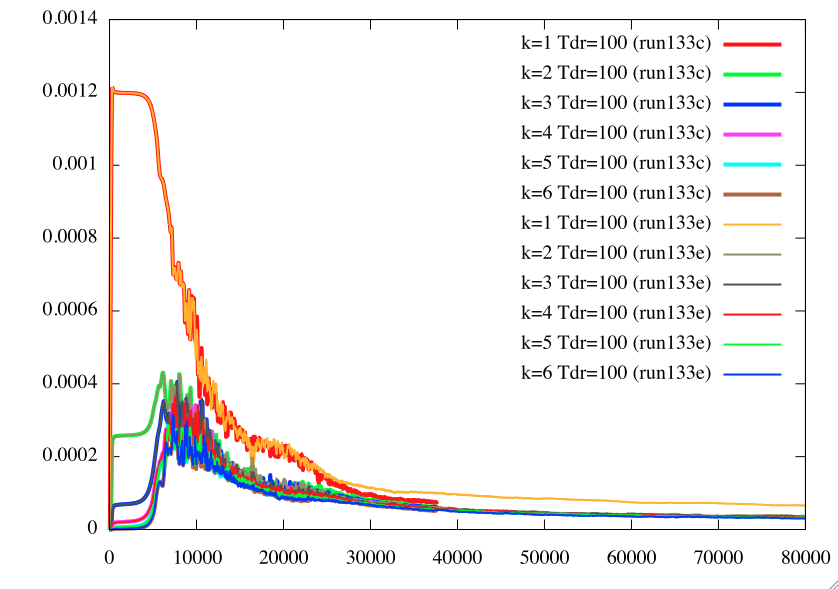}\\
\includegraphics[trim=0.cm 0cm 0cm 0cm, clip=true, width=0.25\linewidth,angle=0,height=5cm]{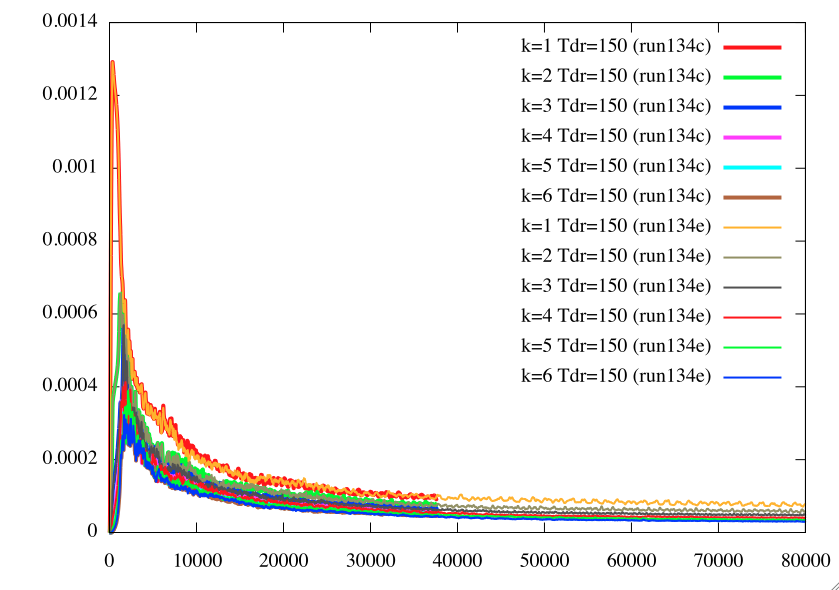}\\
\includegraphics[trim=0.cm 0cm 0cm 0cm, clip=true, width=0.25\linewidth,angle=0,height=5cm]{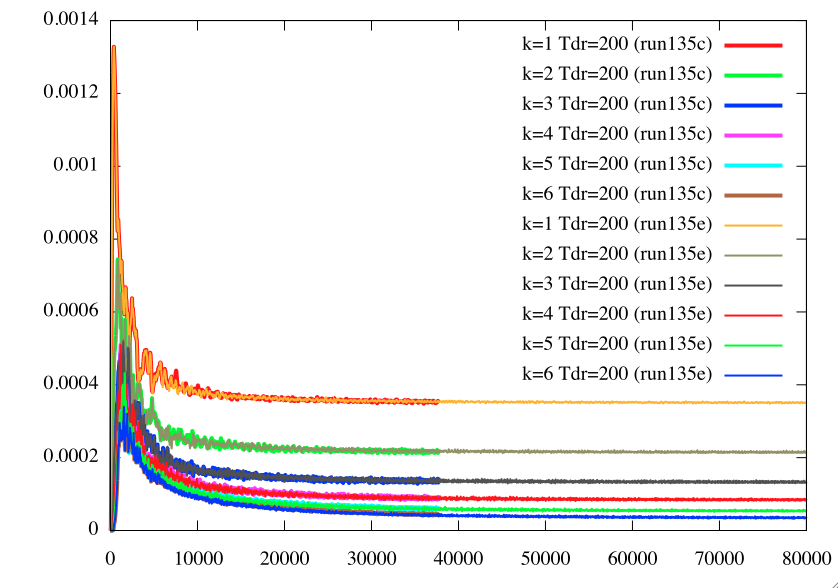}
\end{tabular}
\caption{Time evolution of $\sqrt{\int_\mathbb{R}|\hat{f}^{(k)}(t,v)|^2dv}$ for $k=1,\ 2,\dots, 6$ and for drive time $50,\ 100,\ 150,\ 200$ (run132c-135c and run132e-135e)
Parameters are $\Delta t=0.5$, non-uniform, $6$-th order time scheme and $(N_x,N_v)=(4096,32768)$ for run132c-135c, $(N_x,N_v)=(2048,16384)$ for run132e-135e.}
\label{fig23}
\end{figure}

\begin{figure}[h!]
\begin{tabular}{c}
\includegraphics[trim=0.cm 0cm 0cm 0cm, clip=true, width=0.25\linewidth,angle=0,height=5cm]{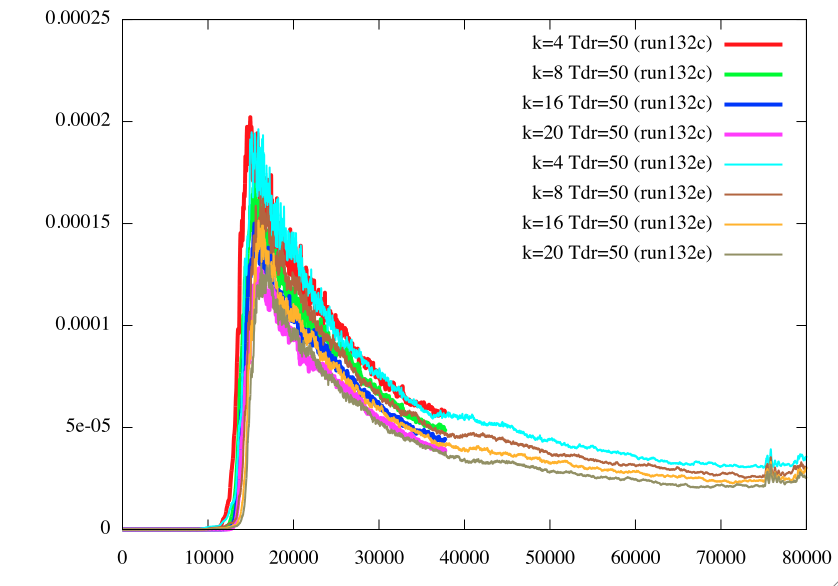}\\
\includegraphics[trim=0.cm 0cm 0cm 0cm, clip=true, width=0.25\linewidth,angle=0,height=5cm]{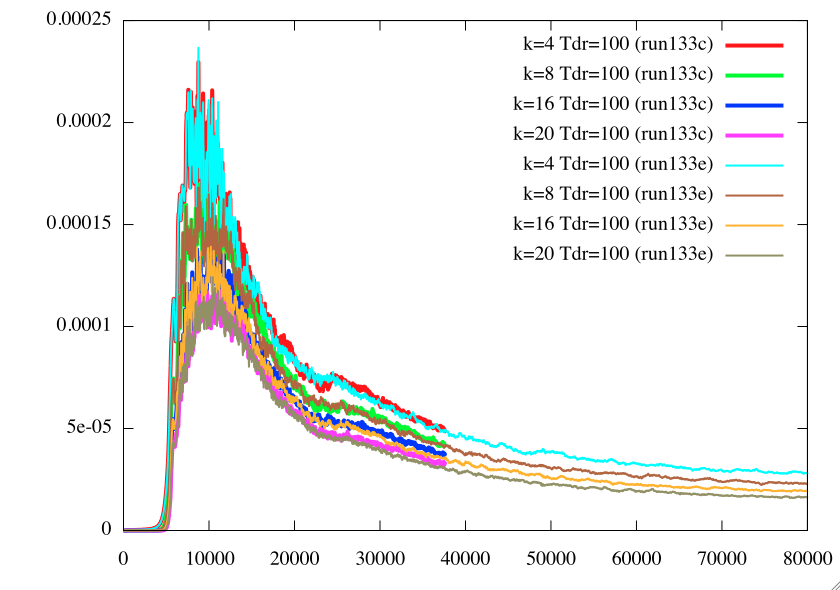}\\
\includegraphics[trim=0.cm 0cm 0cm 0cm, clip=true, width=0.25\linewidth,angle=0,height=5cm]{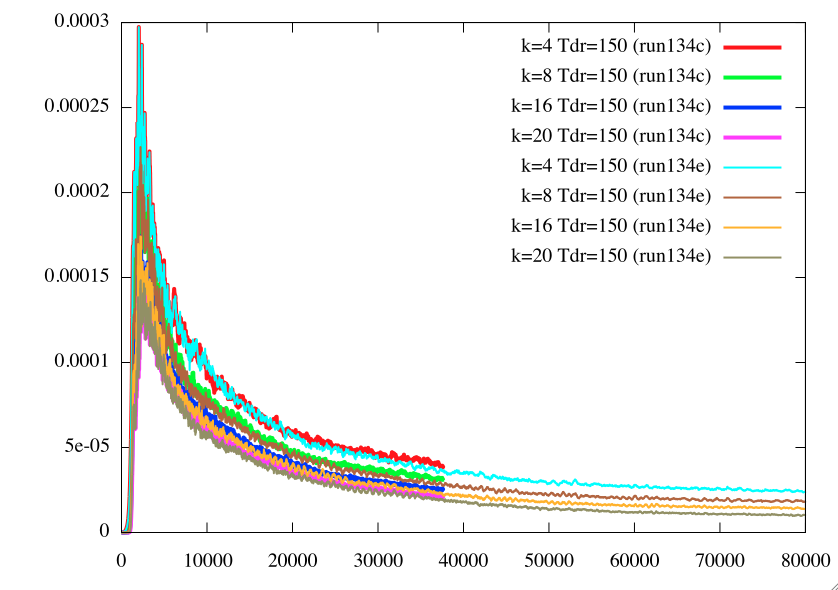}\\
\includegraphics[trim=0.cm 0cm 0cm 0cm, clip=true, width=0.25\linewidth,angle=0,height=5cm]{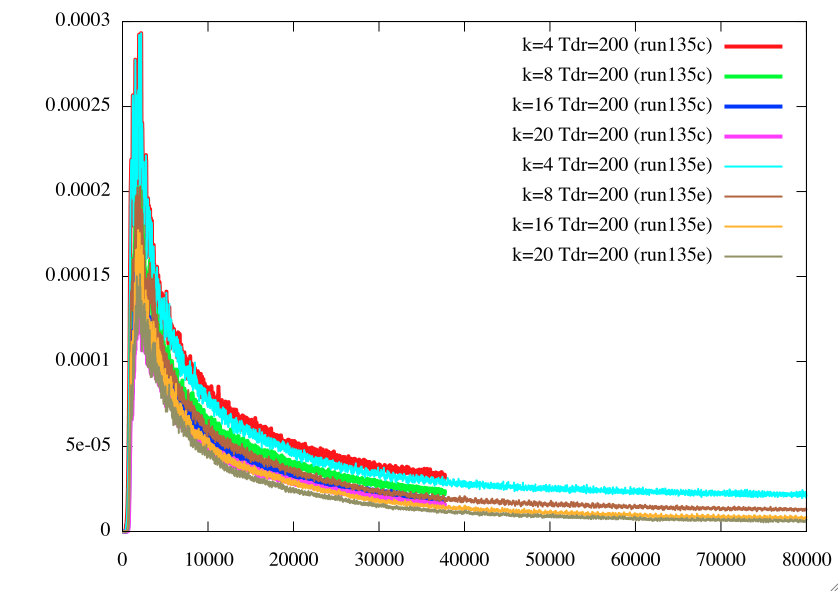}
\end{tabular}
\caption{Time evolution of $\sqrt{\int_\mathbb{R}|\hat{f}^{(k)}(t,v)|^2dv}$ for $k=4,\ 8,\ 16,\ 20$ and for drive time $50,\ 100,\ 150,\ 200$ (run132c-135c and run132e-135e)
Parameters are $\Delta t=0.5$, non-uniform, $6$-th order time scheme and $(N_x,N_v)=(4096,32768)$ for run132c-135c, $(N_x,N_v)=(2048,16384)$ for run132e-135e.}
\label{fig24}
\end{figure}

\begin{figure}[h!]
\begin{tabular}{c}
\includegraphics[trim=5.5cm 6.5cm 1cm 1.7cm, clip=true, width=0.5\linewidth,angle=0,height=7cm]{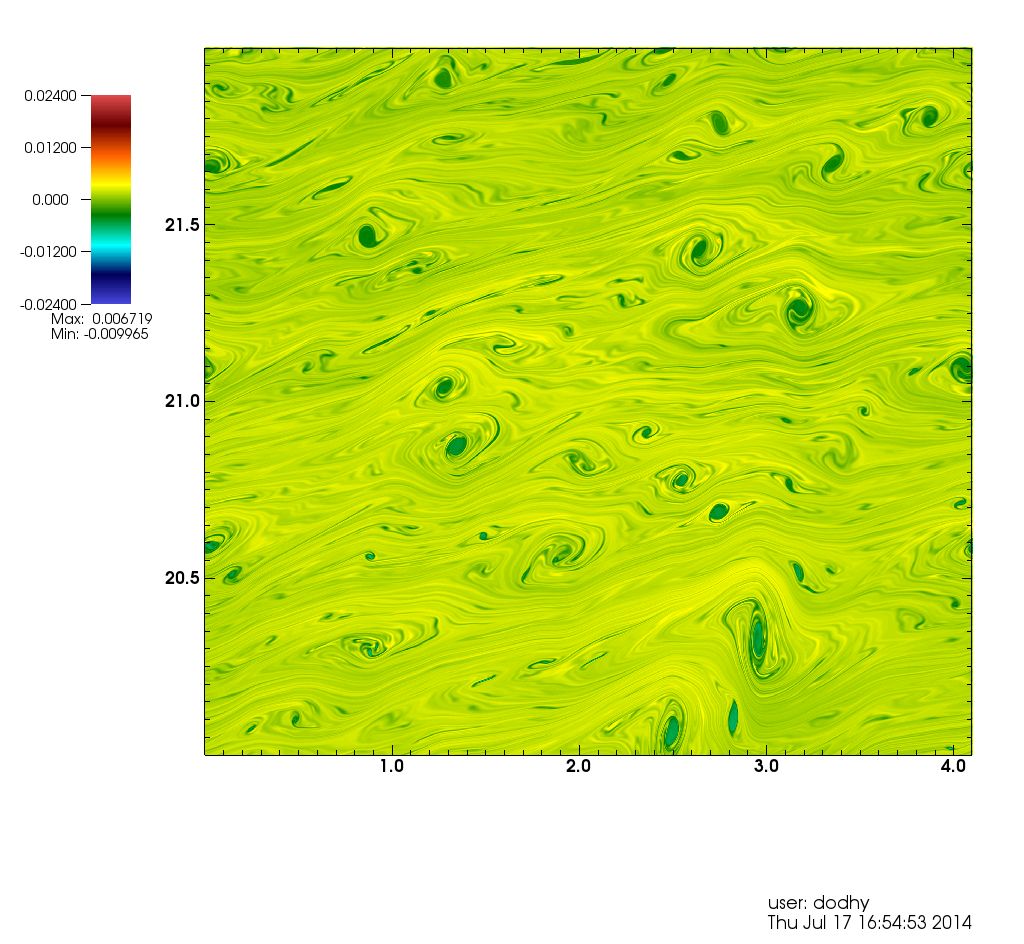}\\
\includegraphics[trim=5.5cm 6.5cm 1cm 1.7cm, clip=true, width=0.5\linewidth,angle=0,height=7cm]{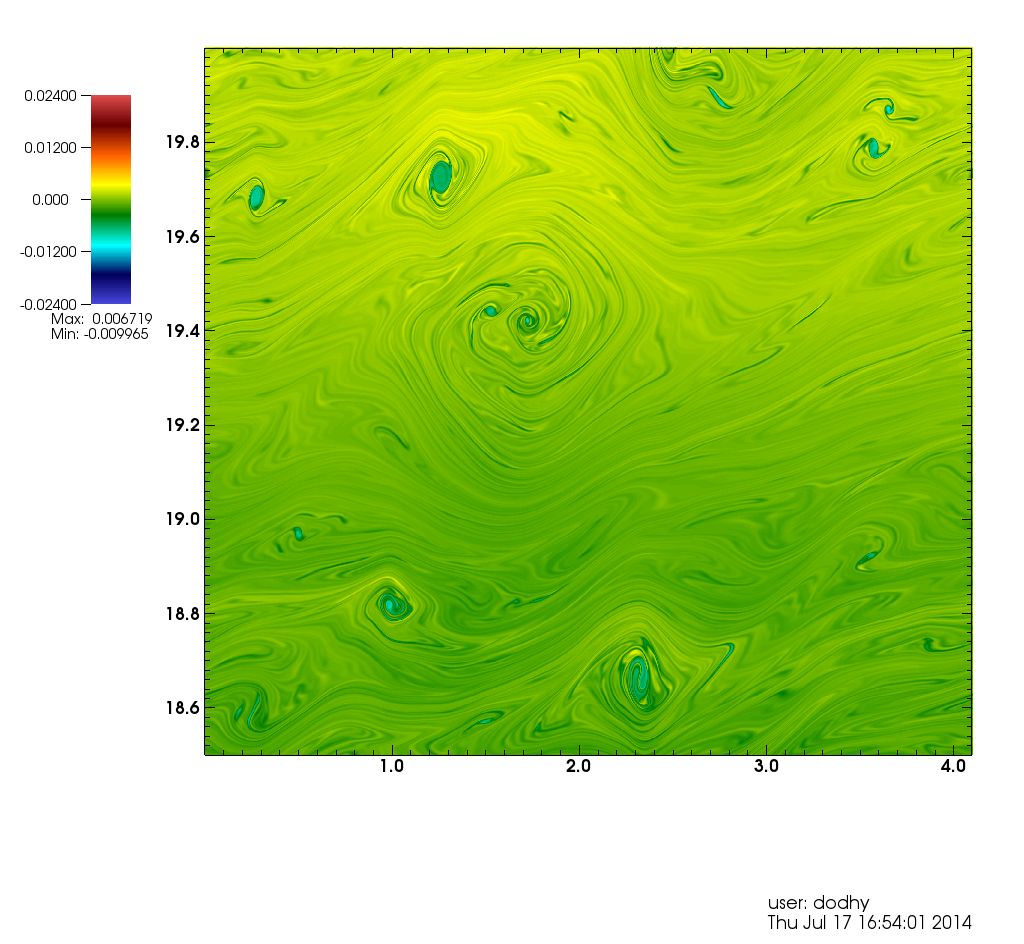}\\
\includegraphics[trim=5.5cm 6.5cm 1cm 1.7cm, clip=true, width=0.5\linewidth,angle=0,height=7cm]{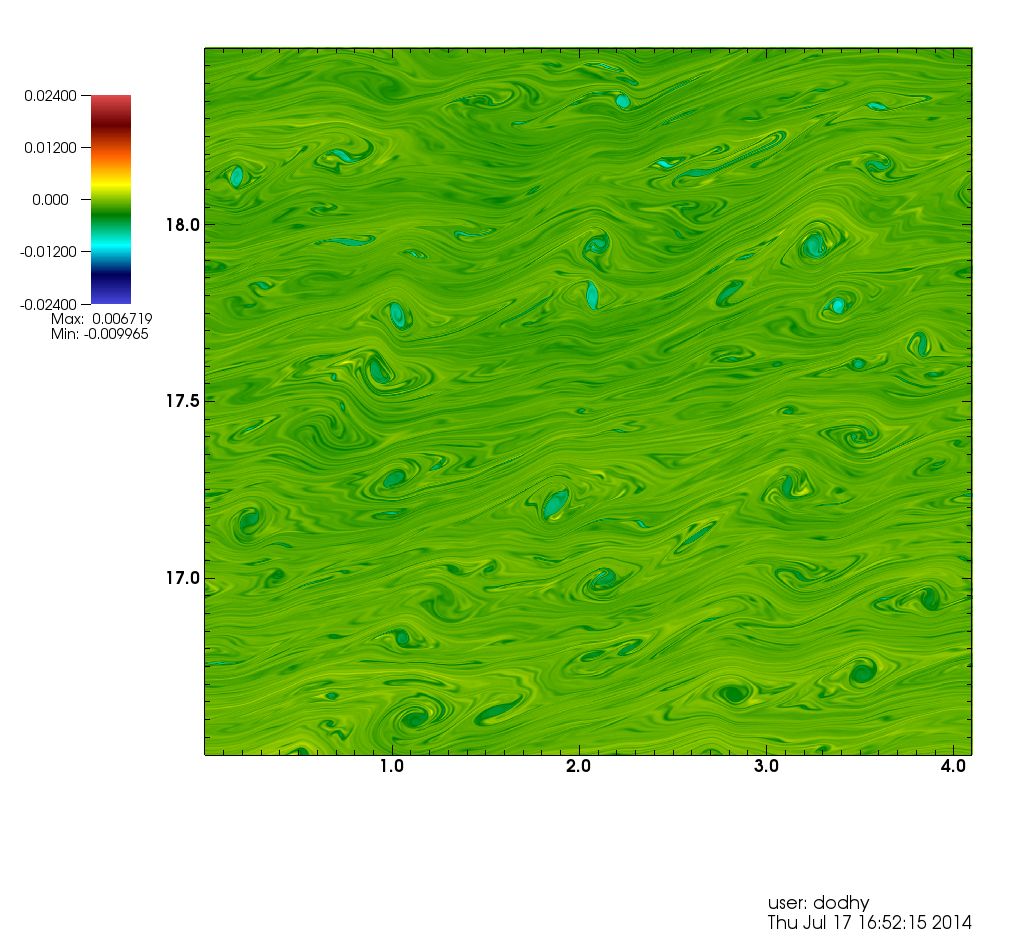}
\end{tabular}
\caption{Details of $\delta f$ distribution function $(f-f_0)(x_i,v_j)$ at time $T=36000$ as a function of $(i,j) \in [0,N_x]\times [16000,22000]$
Parameters are $N_x=4096,\ N_v=32768$,
$\Delta t=0.5$, non-uniform, $6$-th order time scheme, $T_{\rm Dr}=100,\ a_{\rm Dr}=0.00625$ (run133c).}
\label{fig25}
\end{figure}

\begin{figure}[h!]
\begin{tabular}{c}
\includegraphics[trim=5.5cm 6cm 2cm 1cm, clip=true, width=0.5\linewidth,angle=0,height=7.5cm]{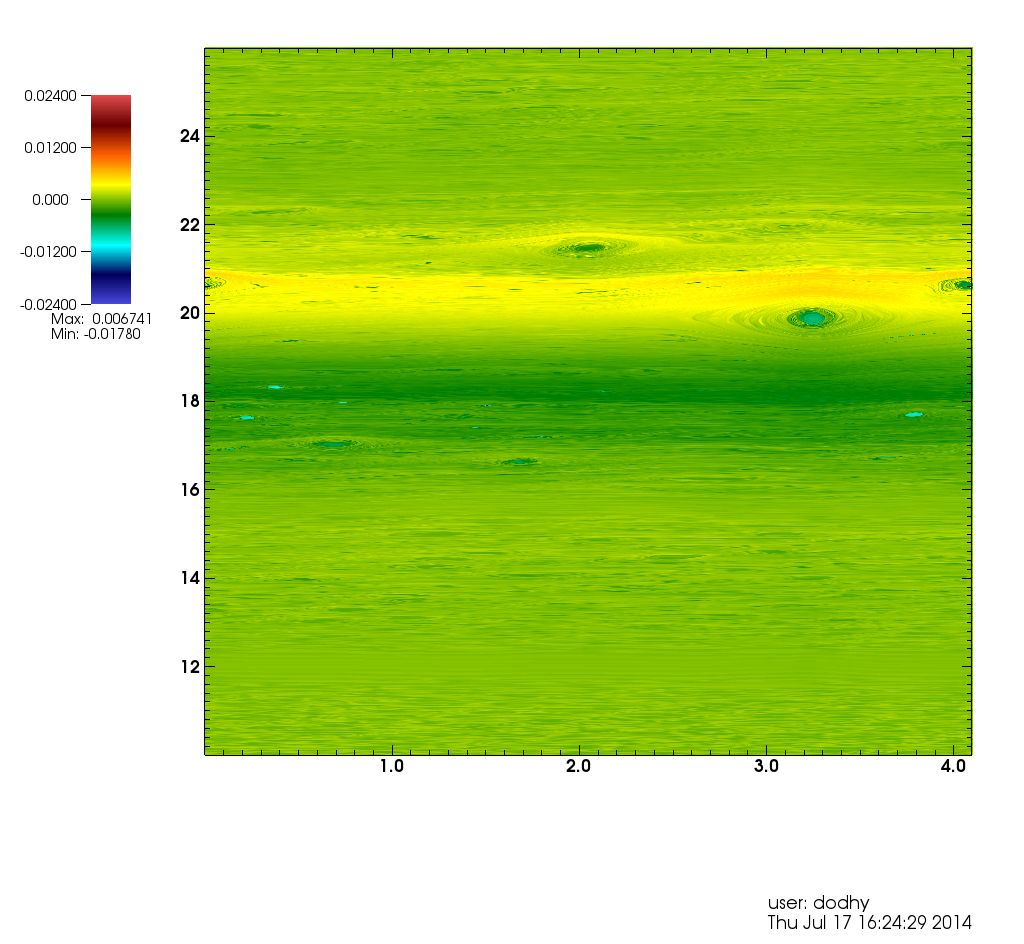}
\end{tabular}
\caption{
$\delta f$ distribution function $(f-f_0)(x_i,v_j)$ at time $T=36000$ as a function of $(i,j) \in [0,N_x]\times [13000,26000]$.
Parameters are $N_x=4096,\ N_v=32768$,
$\Delta t=0.5$, non-uniform, $6$-th order time scheme, $T_{\rm Dr}=150,\ a_{\rm Dr}=0.00625$ (run134c).
}\label{fig26}
\end{figure}

\begin{figure}[h!]
\begin{tabular}{c}
\includegraphics[trim=5.5cm 6.cm 2cm 1.cm, clip=true, width=0.5\linewidth,angle=0,height=7.5cm]{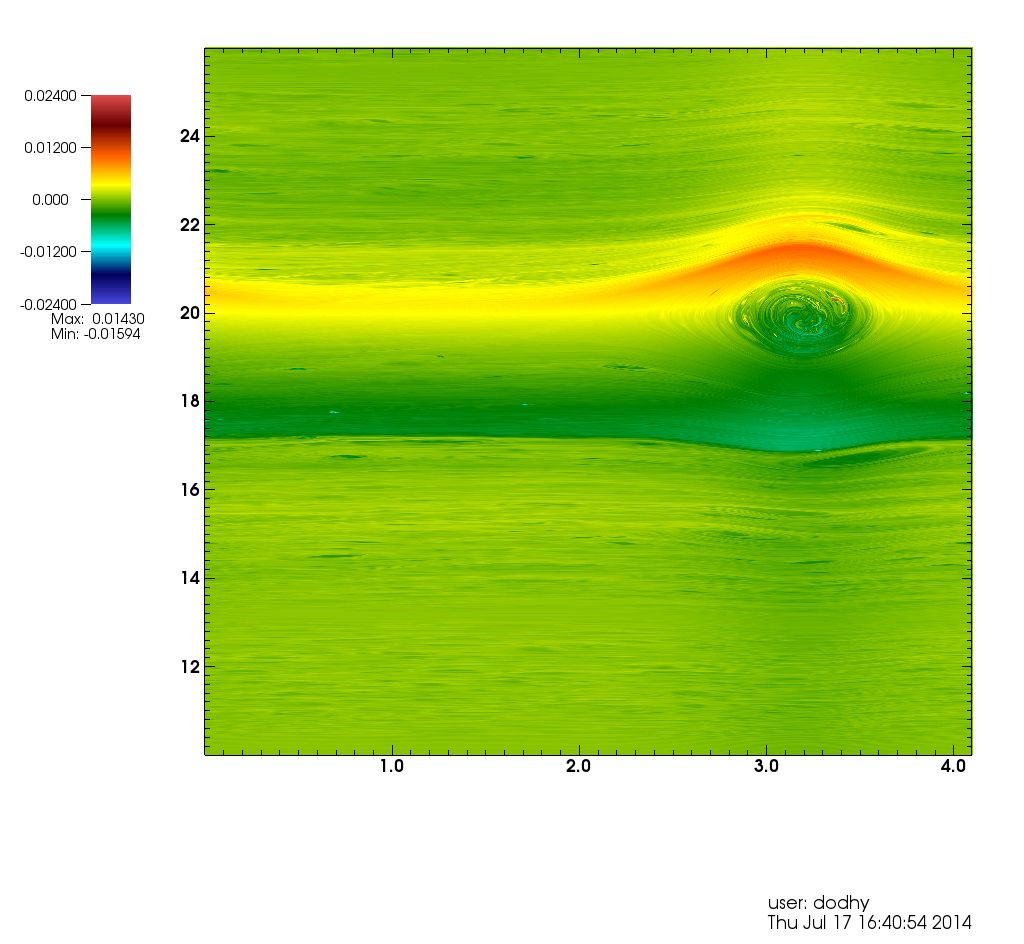}
\end{tabular}
\caption{
$\delta f$ distribution function $(f-f_0)(x_i,v_j)$ at time $T=36000$ as a function of $(i,j) \in [0,N_x]\times [13000,26000]$.
Parameters are $N_x=4096,\ N_v=32768$,
$\Delta t=0.5$, non-uniform, $6$-th order time scheme, $T_{\rm Dr}=200,\ a_{\rm Dr}=0.00625$ (run135c).
%
%
}\label{fig27}
\end{figure}

We study the influence of the drive time, by using the small drive amplitude $a_{\rm Dr}=0.00625$.
Results are given on Figures \ref{Tdrive100} to \ref{fig27}. We can remark that we get formation of the KEEN wave earlier and better formed, by increasing the drive time.
On the contrary, by taking a smaller drive time, phase mixing appears later: see for example on Figure \ref{Tdrive100}, the distribution $\delta f$ function at time $T=5000$ and the $\rho$ harmonics
in semi-logarithmic scale. To complement this diagnostic, we also plot the time evolution of $\sqrt{\int_\mathbb{R}|\hat{f}^{(k)}(t,v)|^2dv}$ for different values of mode $k$, which is another relevant
diagnostic for KEEN waves. We show that with the new diagnostic, we can see a KEEN wave forming when the lower Fourier modes form a definite amplitude ratio pulling away from many modes having the same value with no low mode differentiation. That self-organization in phase space into phase locked low order harmonic states does not happen when there is a proliferation of tiny vortices which have not coalesced, 
as we can see on the delta-f distribution plots (Figures \ref{fig25} to \ref{fig27}). We clearly see the stable KEEN wave for $T_{Dr}=200$, something intermediate with $3$ vortices for $T_{\rm Dr}=150$ and
a forest of vorticlets for $T_{\rm Dr}=100$ at the very late time $T=36000$ where convergence in phase space is no more expected. 
On the other hand, we see how well the RMS quantities are resolved with respect to change of time step (Figures \ref{fig21},\ref{fig22}) or phase space resolution (Figures \ref{fig23},\ref{fig24}).

\section{Conclusion}
The advances made in (i) tackling nonuniform velocity resolution in Vlasov simulations using cubic splines for interpolation, (ii) the use of higher order time splitting schemes, and (iii) efficient parallelization, has allowed accurate long time evolution studies of KEEN waves launched with smaller and smaller drive amplitudes or durations without compromising accuracy. The physics of KEEN wave dynamics is best explored by varying the drive duration, amplitude, wavenumber and frequency. For different choices of ponderomotive drive, $(a_{\rm Dr}, T_{\rm Dr}, \omega_{\rm Dr}$ and $k_{\rm Dr}),$ different birth and spreading of harmonic content physics will be revealed \cite{bedros14a,bedros14b}. Newly discovered features of the dynamics of KEEN waves include the finite partitioning of phase space, unusual particle orbit statistics variations in the topology near separatrices for weak to strong drive, successive harmonic generation instabilities and phase locking, partial mode truncation reconstructions shortcomings and much more \cite{bedros14a,bedros14b}. The key to properly understanding the rich dynamics of Vlasov equations is the adoption of very detailed and physically motivated diagnostics. Some of them were listed above. These need to be developed further in order to automate the analysis of nonlinear self-organized kinetic structures on phase space.   Also, since more than one technique of variable gridding of phase space exists, it is of interest to cross validate these results against, for example, \cite{steiner_vlasovia2013,larson14}. The extension to $2D\times2D$ simulations would then follow. It is also important to extend these techniques to the Vlasov-Maxwell setting and consider stimulated KEEN wave scattering, SKEENS, and not just externally imposed ponderomotive forces. SKEENS and their interaction with SRS will bear fruit in this regard \cite{bedros_ref1}.  It is also imperative to study KEEN-KEEN interactions and KEEN-EPW interactions \cite{bedros14a}. The latter requires fine mesh solutions to encompass the entire region of velocity between the vicinity of the EPW phase velocity down to the vicinity of the well formed KEEN wave phase velocity. This requires of the order of three times as wide a range in which to finely resolved velocity than was used in this paper, where isolated KEEN waves were excited and studied. 



\begin{table}
\begin{tabular}{l}
run, $N_x\times N_v\times \Delta t\times T$, scheme/case, time/proc/eff\\
6, $2^{10}\times 2^{14}\times 1\times 10^4$, $6$-th/U/sm., 2847s/128S/5.1\\
7, $2^{10}\times 2^{15}\times 1\times 10^4$,  $6$-th/U/sm., 5466s/128S/5.3\\
9, $2^{11}\times2^{15}\times0.5\times 5\cdot10^3$, $6$-th/NU/sm., 11866s/128S/4.9\\
10, $2^{11}\times2^{16}\times0.5\times 5\cdot10^3$, $6$-th/NU/sm., 12589s/256S/4.6\\
11, $2^{11}\times 2^{17}\times 0.5\times 5\cdot10^3$, $6$-th/U/sm., 26014s/256S/4.4\\
12, $2^{11}\times2^{18}\times0.5\times 5\cdot10^3$, $6$-th/U/sm., 54800s/256S/4.2\\
13, $2^{11}\times2^{14}\times0.5\times 5\cdot10^3$, $6$-th/NU/sm., 3453s/256S/4.2\\
16a, $2^{11}\times2^{14}\times0.25\times 5\cdot10^3$, $6$-th/NU/sm., 6905s/256S/4.2\\
20, $2^{11}\times2^{14}\times2^{-7}\times 5\cdot10^3$, Str./NU/sm., 72015s/256S/3.5\\
33, $2^{11}\times2^{14}\times 0.5\times 10^4$, $6$-th/NU/sm3, 10898s/256S/5.3\\
45, $2^{11}\times2^{11}\times 0.25\times 10^3$, $6$-th/NU/can, 155s/256M/4.7\\
46, $2^{11}\times2^{12}\times 0.25\times10^3$, $6$-th/NU/can, 284s/256M/5.1\\
47, $2^{11}\times2^{13}\times 0.25\times 10^3$, $6$-th/NU/can, 883s/256M/3.3\\
48, $2^{11}\times2^{14}\times 0.25\times 10^3$, $6$-th/NU/can, 1462s/256M/3.9\\
55, $2^{11}\times2^{17}\times 0.25\times 10^3$, $6$-th/NU/can, 11594s/256M/4.0\\
59, $2^{13}\times 2^{14}\times 1\times 10^3$, $6$-th/NU/can, 1423s/256M/4.1\\
68, $2^{14}\times2^{15}\times0.25\times 10^3$, $6$-th/NU/can, 24509s/256M/3.8\\
72, $2^{11}\times2^{12}\times 0.125\times 10^3$, $6$-th/NU/can, 567s/256M/5.1\\
73, $2^{11}\times 2^{12}\times 0.25\times 10^3$, $6$-th/NU/can, 285s/256M/5.1\\
74, $2^{12}\times 2^{13}\times 0.125\times 10^3$, $6$-th/NU/can, 2766s/256M/4.2\\
75, $2^{12}\times 2^{13}\times 0.25 \times 10^3$, $6$-th/NU/can, 1384s/256M/4.2\\
76, $2^{13}\times 2^{14}\times 0.125\times 10^3$, $6$-th/NU/can, 11255s/256M/4.1\\
77, $2^{13}\times2^{14}\times0.25\times 10^3$, $6$-th/NU/can, 5624s/256M/4.1\\
78, $2^{10}\times 2^{11}\times 0.125\times 10^3$, $6$-th/NU/can, 186s/256M/3.9\\
79, $2^{10}\times 2^{11}\times 0.25 \times 10^3$, $6$-th/NU/can, 94.5s/256M/3.8\\
80, $2^{9}\times 2^{10}\times 0.125 \times 10^3$, $6$-th/NU/can, 112s/256M/1.6\\
81, $2^{9}\times 2^{10}\times 0.25\times 10^3$, $6$-th/NU/can, 55.9s/256M/1.6\\
84, $2^{13}\times 2^{14}\times 0.025\times 10^3$, Str./NU/can, 19018s/256M/3.3\\
85, $2^{13}\times 2^{14}\times 0.05\times 10^3$, Str./NU/can, 9518s/256M/3.3\\
86, $2^{13}\times 2^{14}\times 0.25\times 10^3$, $6$-th/NU/can, 5567s/256M/4.1\\
87, $2^{13}\times 2^{15}\times 0.25\times 10^3$, $6$-th/NU/can, 11469s/256M/4.0\\
88, $2^9\times 2^{12}\times 0.05\times 10^3$, Str./U/can, 169s/256M/2.9\\
89, $2^{13}\times2^{14}\times0.0125\times 10^3$, Str./NU/can, 35545s/256M/3.5\\
90, $2^{13}\times2^{16}\times0.25\times 10^3$, $6$-th/U/can, 33831s/256M/2.7\\
102, $2^{13}\times2^{16}\times0.125\times 10^3$, $6$-th/NU/can, 48852s/256M/3.8\\
104, $2^{13}\times 2^{14}\times 0.5 \times 10^3$, $6$-th/NU/can, 2828s/256M/4.1\\
105, $2^{13}\times 2^{14}\times 0.1\cdot2^{-4} \times 10^3$, Str./NU/can, 70164s/256M/3.6\\
106b, $2^{13}\times 2^{14}\times 0.1\cdot2^{-4}\times 10^3$, $6$-th/NU/can, 46510s/1024S/5.0\\
107b, $2^{13}\times 2^{14}\times 0.1\cdot2^{-5}\times 923$, $6$-th/NU/can, 24h/1024S/4.9\\
108c, $2^{13}\times 2^{14}\times 0.1\cdot2^{-6}\times 683$, $6$-th/NU/can, 24h/2048S/3.6\\
132c, $2^{12}\times 2^{15}\times 0.5\times 37759.5$, $6$-th/NU/sm0.5, 24h/256M/5.0\\
133c, $2^{12}\times 2^{15}\times 0.5\times 37617.5$, $6$-th/NU/sm1, 24h/256M/5.0\\
134c, $2^{12}\times 2^{15}\times 0.5\times 37627$, $6$-th/NU/sm1.5, 24h/256M/5.0\\
135c, $2^{12}\times 2^{15}\times 0.5\times 37726$, $6$-th/NU/sm, 24h/256M/5.0\\
132d, $2^{12}\times 2^{15}\times 0.25\times 18875.75$, $6$-th/NU/sm0.5, 24h/256M/5.0\\
133d, $2^{12}\times 2^{15}\times 0.25\times 18911.5$, $6$-th/NU/sm1, 24h/256M/5.0\\
134d, $2^{12}\times 2^{15}\times 0.25\times 18846.75$, $6$-th/NU/sm1.5, 24h/256M/5.0\\
135d, $2^{12}\times 2^{15}\times 0.25\times 18828.5$, $6$-th/NU/sm, 24h/256M/5.0\\
132e, $2^{11}\times 2^{14}\times 0.5\times 8\cdot10^4$, $6$-th/NU/sm0.5, 51015s/256M/4.5\\
133e, $2^{11}\times 2^{14}\times 0.5\times 8\cdot10^4$, $6$-th/NU/sm1, 51071s/256M/4.5\\
134e, $2^{11}\times 2^{14}\times 0.5\times 8\cdot10^4$, $6$-th/NU/sm1.5, 51132s/256M/4.5\\
135e, $2^{11}\times 2^{14}\times 0.5\times 8\cdot10^4$, $6$-th/NU/sm, 51099s/256M/4.5\\
216, $2^{11}\times 2^{14}\times 0.5\times 10^4$, $6$-th/NU/sm1, 6300s/256H/4.6\\ 
294, $2^{13}\times 2^{14}\times 0.025\times 10^3$, $6$-th/NU/can, 43920s/256H/5.3
\end{tabular}
\caption{Parameters for the different runs.
The first column, such as $6$-th, refers to the order of the time splitting scheme.Str. is for Strang splitting.
U is for uniform/ NU for non-uniform velocity grids. sm. is for small drive cases ($a_{\rm Dr}=0.00625$ and $T_{\rm Dr}=200$); sm1 is for $a_{\rm Dr}=0.00625$ and $T_{\rm Dr}=100$;
sm3 is for $a_{\rm Dr}=0.00625$ and $T_{\rm Dr}=300$; sm1.5 is for $a_{\rm Dr}=0.00625$ and $T_{\rm Dr}=150$ and so on;
can is for the canonical drive case ($a_{\rm Dr}=0.2$ and $T_{\rm Dr}=100$).
S is for selavlas allocation on Helios; M is for mic\_eu allocation on Helios; H is for Hydra; eff is for efficiency, that is
$(s\cdot N_x\cdot N_v\cdot T/\Delta t)/(\textrm{time (in s.)}\cdot10^{6}\cdot \textrm{proc})$, where $s=3$ for Strang and $s=11$ for the $6$-th order scheme.}
\label{table1}
\end{table}

\clearpage

\clearpage

%
%

\clearpage

\clearpage

\newpage


\clearpage

\clearpage

\begin{figure}[h!]
\begin{center}
\begin{tabular}{ccc}
\includegraphics[trim=0.cm 0cm 0cm 0cm, clip=true, width=0.15\linewidth,angle=0,height=3.2cm]{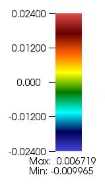} &
\includegraphics[trim=0.cm 0cm 0cm 0cm, clip=true, width=0.15\linewidth,angle=0,height=3.2cm]{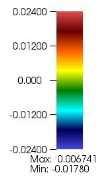} &
\includegraphics[trim=0.cm 0cm 0cm 0cm, clip=true, width=0.15\linewidth,angle=0,height=3.2cm]{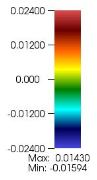}
\end{tabular}
\end{center}
\caption{Legends for Figures \ref{fig25},\ref{fig26},\ref{fig27}.}
\end{figure}

\end{document}